\documentclass{nm}
\usepackage{graphicx}
\usepackage{subfigure}
\usepackage{algorithm}
\usepackage{algorithmic}
\usepackage{mathtools}

\newcommand{\mc}{\mathcal}
\newtheorem{dfn}{Definition}[section]
\newcommand{\hpsngAbs}{\ensuremath{\!\!\!\!\!\!\!\!\! = \ \ \ }}
\newcommand{\hpsngabs}{\ensuremath{\!\!\!\!\!\! = \ \ }}
\newcommand{\hiint}{\ensuremath{\iint \!\!\!\!\!\!\! = \ \ \!\!\!\!\!\!\!\!\!\!\!\!\! = \ \ }}
\newcommand{\hint}{\ensuremath{\int \!\!\!\!\!\!\! = \ \ }}
\DeclareMathOperator\supp{supp}
\newcommand{\mathbbm}[1]{\text{\usefont{U}{bbm}{m}{n}#1}} % from mathbbm.sty

\begin{document}

\markboth{X. Feng, Y. Qian and W. Shen}{MC-Nonlocal-PINNs}
\title{MC-Nonlocal-PINNs: handling nonlocal operators in PINNs via Monte Carlo sampling}

\author[X. Feng, Y. Qian and W. Shen]{Xiaodong Feng\affil{1}\comma\corrauth, Yue Qian\affil{1} and  Wanfang Shen \affil{2}}

\address{\affilnum{1}\ Institute of Computational Mathematics and Scientific/Engineering
	Computing, Academy of Mathematics and Systems Science, Chinese Academy
	of Sciences, Beijing, China. \\
	\affilnum{2}\ Shandong Key Laboratory of Blockchain Finance, Shandong University of Finance and Economics,
	Jinan 250014, China.\\}

\emails{{\tt xdfeng@lsec.cc.ac.cn} (X. Feng), {\tt qianyue2021@lsec.cc.ac.cn} (Y. Qian), {\tt wfshen@sdufe.edu.cn} (W. Shen)}
	%%%%% Begin Abstract %%%%%%%%%%%
	\begin{abstract}
		We propose, Monte Carlo Nonlocal physics-informed neural networks (MC-Nonlocal-PINNs), which is a generalization of MC-fPINNs in \cite{guo2022monte}, for solving general nonlocal models such as integral equations and nonlocal PDEs. Similar as in MC-fPINNs, our MC-Nonlocal-PINNs handle the nonlocal operators in a Monte Carlo way, resulting in a very stable approach for high dimensional problems. We present a variety of test problems, including high dimensional Volterra type integral equations, hypersingular integral equations and nonlocal PDEs, to demonstrate the effectiveness of our approach.
	\end{abstract}
	%%%%% end %%%%%%%%%%%
	
	%%%%% Keywords %%%%%%%%%%%
	\keywords{Nonlocal models,  PINNs, Monte Carlo sampling, deep neural networks.}

	%%%% AMS subject classifications %%%%
	\ams{65C05, 65D30, 65R20}
	
	%%%% maketitle %%%%%
	\maketitle
	
	%%%% Start %%%%%%

	\section{Introduction}
	
	Deep neural networks have gained a growing interest in recent years with a wide variety of methods ranging from computer vision and natural language processing to simulations of physical systems \cite{weinan2020machine, lu2021deepxde, yu2018deep, sirignano2018dgm}. A representative example is physics-informed neural networks(PINNs) \cite{raissi2019physics}, whose central idea is to incorporate governing laws of physical systems into the training loss function and recast the original problem into an optimization problem. PINNs have demonstrated remarkable success in applications including fluid mechanics \cite{raissi2020hidden,brunton2020machine}, high dimensional PDEs (with applications in computational finance) \cite{han2018solving,zang2020weak, huang2022augmented}, uncertainty quantification \cite{yang2021b,qin2021deep,zhang2018deep,iten2020discovering,meng2020composite, guo2022normalizing}, to name just a few.

	For PDE models with classic (integer) derivatives, PINNs adopt automatic differentiation to solve PDEs by penalizing the PDE in the loss function at a
	random set of points in the domain of interest.  However, for PDE models involving nonlocal operators, one can no longer use the automatic differentiation to handle the operators due to the nonlocal property. To overcome this issue, fPINNs \cite{pang2019fpinns} was developed for solving space-time fractional advection-diffusion equations. The main idea in \cite{pang2019fpinns} is to introduce a classic discretization technique to handle the fractional operator. However, this is not a good choice for high dimension problems since the curse of dimensionality. Similar idea has been used to handle more general non-local operators in \cite{pang2020npinns}, while the approach again can not be used for high dimensional cases. We also mention the work \cite{yuan2022pinn}, where the so called A-PINN was proposed to handle some special types of integral equations.

	More recently, the MC-fPINNs approach was proposed in \cite{guo2022monte} to handle fractional PDEs, where the fractional operators are handled in a Monte Carlo way, resulting in a very stable approach for high dimensional problems. Take the fractional Laplacian equation as an example:
	\begin{equation}
		(-\Delta)^{\alpha/2}u(x)=C_{d,\alpha} \mbox{ P.V.} \int_{\mathbb{R}^d}\frac{u(x)-u(y)}{\Vert x-y\Vert _2^{d+\alpha}}dy,\quad 0<\alpha<2,
		\label{laplacian}
	\end{equation}
	where \mbox{P.V.} denotes the principle value of the integral and $C_{d,\alpha}$ is a constant depending on $\alpha$ and $d$. One can divide the integral into the following two parts:
	\begin{equation}\label{laplacian_two_parts}
		(-\Delta)^{\alpha/2} u(x)=C_{d,\alpha}\left(\int _{y\in B_{r_0}(x)}\frac{u(x)-u(y)}{\Vert x-y\Vert_2^{d+\alpha}}dy + \int_{y\notin B_{r_0}(x)}\frac{u(x)-u(y)}{\Vert x-y\Vert_2^{d+\alpha}}dy \right).
	\end{equation}
	It is shown that the fractional Laplacian of $u(x)$ can be calculated via the following approximation:
	\begin{equation}
		\begin{aligned}
			(-\Delta)^{\alpha/2}u(x)=&C_{d,\alpha} \frac{\left|S^{d-1}\right|r_0^{2-\alpha}}{2(2-\alpha)} \mathbb{E}_{\xi,r_I\sim f_I(r)} \left[\frac{2u(x)-u(x-r_I\xi)-u(x+r_I\xi)}{r_I^2}\right] \\
			&+ C_{d,\alpha} \frac{\left|S^{d-1}\right|r_0^{-\alpha}}{2\alpha}\mathbb{E}_{\xi,r_O\sim  f_O(r)}\left[2u(x)-u(x-r_O\xi) -u(x+r_O \xi)\right],
		\end{aligned}
		\label{laplacian_integral}
	\end{equation}
	where $\left|S^{d-1}\right|$ denotes the surface area of $S^{d-1}$, $\xi$ is uniformly distributed on the sphere $S^{d-1}$, and $r_I,r_O$ can be quickly sampled via
	\begin{equation}
		r_I/r_0 \sim \mbox{ Beta}(2-\alpha, 1),\quad r_0/r_O \sim \mbox{ Beta}(\alpha,1).
	\end{equation}
	More precisely, one can resort to the classic Monte Carlo sampling to handle the fractional Laplacian (see in \cite{guo2022monte} for more details).	
	
	The main aim of this work is to extend the idea in \cite{guo2022monte} to more general nonlocal operators. Our new contributions are summarized as follows:
	\begin{itemize}
		\item We generalize MC-fPINNs to MC-Nonlocal-PINNs, which can handle more general nonlocal models such as Volterra type integral equations with either bounded or singular kernels, hypersingular integral equations, and nonlocal PDEs with various kernels.
		\item We present several high dimensional examples to show the effectiveness of the MC-Nonlocal-PINNs approach.
	\end{itemize}
	
	The remainder of this paper is structured as follows. In Section 2, we present some preliminaries. In Section 3, we present our MC-Nonlocal-PINNs approach for solving general nonlocal problems. This is followed by numerical tests in Section 4. Finally, we give some concluding remarks in Section 5.

	\section{Preliminaries}
	
	In this section, we present some preliminaries.
	\subsection{Physics-informed neural network}
	In this section, we first give a brief review on physics informed neural networks (PINNs) \cite{raissi2019physics, pang2019fpinns}.  To this end, we consider the following PDE
	\begin{equation}
		\label{eqn_general}
		\left\{
		\begin{aligned}
			\mathcal{L}_x[u](x)=f(x),\quad x\in\Omega,\\
			u(x) = g(x),\quad x\in \partial \Omega.
		\end{aligned}
		\right.
	\end{equation}
	Here $\mathcal{L}_x$ is a differential operator, and $\Omega$ is a domain of interest. A deep neural network (DNN) is a sequence alternative composition of linear functions and nonlinear activation function. The PINNs approach uses the output of DNN, $u_{NN}(x;\theta)$, to approximate the solution of equation $u(x)$ and calculate the differential operator via automatic differentiation. Here $\theta$ is a collection of the all learnable parameters in the DNN. Specifically, define the PDE residual as
	\begin{equation}
		r(x;\theta) = \mathcal{L}_xu_{NN}(x;\theta) - f(x),
	\end{equation}
	then $\theta$ can be learned by minimizing the following composite loss function
	\begin{equation}\label{loss}
		\mathcal{L}(\theta) = w_r\cdot \mathcal{L}_r(\theta) +w_b\cdot  \mathcal{L}_b(\theta),
	\end{equation}
	where
	\begin{equation}
		\begin{aligned}
			\mathcal{L}_r(\theta) = \frac{1}{2}\sum_{i=1}^{N_r}\big\vert r(x_r^i;\theta)\big\vert ^2 \quad \mbox{and}\quad \mathcal{L}_b(\theta) = \frac{1}{2}\sum_{i=1}^{N_b}\big\vert u_{NN}(x_b^i;\theta)-g(x_b^i)\big\vert ^2,
		\end{aligned}
	\end{equation}
	$\{w_r,w_b\}$ are weights and $\{x_r^i\}_{i=1}^{N_r}$, $\{x_b^i\}_{i=1}^{N_b}$ denote the training data.
	%{\color{red}If ... local --- AD ... otherwise \cite{fPINN} }
	
	If $\mathcal{L}_x$ in Eq.(\ref{eqn_general}) is a local differential operator, we can easily compute its value via automatic differentiation; otherwise we must first deal with the operator $\mathcal{L}_x$ due to the nonlocal property. Here we remark that the fPINNs \cite{pang2019fpinns} was proposed to solve fractional advection-diffusion equations.
	\subsection{MC-fPINNs}
	Notice that the nonlocal properties bring difficulties for solving PDEs via automatic differentiation, here we briefly introduce the MC-fPINNs in \cite{guo2022monte}, which handle fractional Laplacian operator by a directly Monte Carlo sampling. Let $\Omega\in\mathbb{R}^d$ be a spatial domain, and we denote by $x\in\Omega$ the spatial variable. Consider the standard fractional Laplacian equation
	\begin{equation}
		\begin{aligned}
			& 	(-\Delta)^{\alpha/2}u(x) = f(x), \quad x\in \Omega,\\
			& u(x) = g(x),\quad x\in \mathbb{R}^d \backslash \Omega,
		\end{aligned}
		\label{fractional_pde}
	\end{equation}
	where $0<\alpha<2$, the fractional Laplacian operator $\Delta^{\alpha/2}u(x)$ is defined by (\ref{laplacian}) and $f(x), g(x)$ are given functions. The MC-fPINNs solve the above PDE problems via constructing a DNN model $u_{NN}(x;\theta)$, parametrized by $\theta$, to approximate the solution $u(x)$. Specifically, notice that the fractional operator can be represented as an integral formulation (\ref{laplacian_integral}), hence the fractional operator of the $u_{NN}$ can be calculated as follows:
	\begin{equation}
		\begin{aligned}
			(-\Delta)^{\alpha/2}u_{NN}(x;\theta)=&C_{d,\alpha} \frac{\left|S^{d-1}\right|r_0^{2-\alpha}}{2(2-\alpha)} \mathbb{E}_{\xi,r_I\sim f_I(r)} \left[\frac{2u_{NN}(x;\theta)-u_{NN}(x-r_I\xi;\theta)-u_{NN}(x+r_I\xi;\theta)}{r_I^2}\right] \\
			&+ C_{d,\alpha} \frac{\left|S^{d-1}\right|r_0^{-\alpha}}{2\alpha}\mathbb{E}_{\xi,r_O\sim  f_O(r)}\left[2u_{NN}(x;\theta)-u_{NN}(x-r_O\xi;\theta) -u_{NN}(x+r_O \xi;\theta)\right],
		\end{aligned}
		\label{laplacian_integral_surrogate}
	\end{equation}
	where $C_{d,\alpha}, S^{d-1}, r_0, r_I, r_O, f_I, f_O$ are as prescribed before. Then the residual loss of the PDE (\ref{fractional_pde}) can be written as
	\begin{equation}
		\begin{aligned}
			\mathcal{L}_r(\theta) = \sum_{i=1}^{N_r} \big|(-\Delta)^{\alpha/2}u_{NN}(x_r^i;\theta) - f(x_r^i)\big|^2,
		\end{aligned}	
	\end{equation}
	where $\{x_r^i\}_{i=1}^{N_r}$ denotes training data and $(-\Delta)^{\alpha/2}u_{NN}(x;\theta)$ is computed via (\ref{laplacian_integral_surrogate}).

\section{MC-Nonlocal-PINNs} \label{section:methodology}
In this section, we present our MC-Nonlocal-PINNs approach, which is a generalization of the original MC-fPINNs. We shall mainly consider three typical types of models: Volterra equations, hypersingular equations and nonlocal PDEs.

\subsection{Volterra type equations}
Integral and integro-differential equations, in which the unknown function appears inside an integral sign, have been widely employed in different fields, such as population growth \cite{apreutesei2009travelling}, acoustic scattering \cite{colton2013integral, antoine2021introduction}, mechanics and plasma physics \cite{meleshko2010symmetries}. In this section, we consider a framework for volterra integral and integro-differential equations. A standard \textit{integral equation} has the form of
\begin{equation}
	u(x) = f(x)+\lambda \cdot \int_{h_1(x)}^{h_2(x)}K(x,s)u(s)ds,
	\label{integral_equation}
\end{equation}
where $h_1(x)$ and $h_2(x)$ are the limits of integration, $\lambda$ is a constant parameter, $K(x,s)$ is called the \textit{kernel} of the integral equation. Here functions $f(x)$ and $K(x,s)$ are given in advance and $u(x)$ is the unknown quantity. Without loss of generality, we assume that
$K(x,s)\geq 0.$

An \textit{integro-differential equation} contains an additional derivative operator compared with the original integral equation

\begin{equation}
	\mathcal{N}_x[u](x) = f(x) + \int_{h_1(x)}^{h_2(x)} K(x,s)u(s)ds,
	\label{integro_differential_equation}
\end{equation}
where $\mathcal{N}_x$ denotes a differential operator with respect to $x$ and others are as prescribed before.

To ease the discussion we use notation IDEs to express integral and integro-differential equations. The limits of integration are used to characterize IDEs. When $h_1(x)$ and $h_2(x)$ are fixed(independent of $x$), the form of Eq. (\ref{integro_differential_equation}) is called Fredholm equation; when at least one of $h_1(x)$ and $h_2(x)$ is variable, the form of Eq. (\ref{integro_differential_equation}) is called Volterra equation. If the equation contains nonlinear functions of $u(x)$, such as $\sin(u), e^{u}, \ln(1+u)$, the IDEs are called \textit{nonlinear}. In this work, we focus on forward (non-)linear IDEs, including Fredholm and Volterra types.

\subsection{Hypersingular integral equations}
Many physical problems can be modeled by boundary integral equations with Hadamard-type hypersingular kernels, such as acoustic and solid mechanics \cite{de2005hypersingular,wu2008toeplitz,li2010newton}. The concept of hypersingular integrals was introduced by Hadamard which is defined by the limit of an expansion, ignoring those diverging terms.

\begin{dfn}\label{def1}
	Assume that $u$ is a function defined on $(0,\beta)$ and that there exists the following expansion:
	\begin{equation}
		u(\epsilon) = \sum_{n=0}^N\sum_{m=0}^{M_n}u_{nm}\epsilon^{\tau_n}\log^m\epsilon +U(\epsilon),
	\end{equation}
	where
	\begin{equation}
		\tau_N\leq \tau_{N-1}\leq \cdots\leq \tau_1\leq \tau_0=0
	\end{equation}
	and $u_{j0}=0$ if $\tau_j=0(0\leq j\leq N)$. If $\lim_{\epsilon\to 0}U(\epsilon)$ exists then the finite-part limit of $f(\epsilon)$ as $\epsilon\to 0$ is defined by
	\begin{equation}
		\mbox{F.P. }\lim\limits_{\epsilon\to 0}u(\epsilon) = \lim\limits_{\epsilon\to 0}U(\epsilon).
	\end{equation}
\end{dfn}

One of the major problems arising numerical methods is how to evaluate the following hypersingular integral efficiently
\begin{equation}
	\mathcal{I}(u,s) \coloneqq \int_a^b \hpsngAbs \frac{u(x)}{(x-s)^2}\,dx, \quad s\in(a,b),
	\label{hypersingular_def}
\end{equation}
where $\int\hpsngabs$ denotes a hypersingular integral and $s$ is the singular point. Using the Definition \ref{def1}, we have \begin{equation}
	\int_a^b\hpsngAbs \frac{u(x)}{(x-s)^2}\mathrm{d}x=\mbox{F.P. }\lim\limits_{\epsilon\to 0}\left\{
	\int_a^{s-\epsilon}\frac{u(x)}{(x-s)^2}dx +\int_{s+\epsilon}^{b}\frac{u(x)}{(x-s)^2}dx
	\right\}.
	\label{hypersingular_limit_def}
\end{equation}
Similarly, two dimensional hypersingular integrals can be derived. Without loss of generality, we consider the two dimensional region $\Omega$ with a boundary described by the equation $R=R(\nu),0\leq \nu \leq 2\pi$, with origin point $(0,0)$ of $\Omega$. We consider the following hypersingular integral
\begin{equation}
	{\hiint}_{\Omega}\frac{u(x_1,x_2)}{r^3}dx_1dx_2 = \int_0^{2\pi}\left[{\hint}_0^{R(\nu)}\frac{u(r\cos\nu, r\sin\nu)}{r^2}dr\right]d\nu,
	\label{two_dim_hyper_singular}
\end{equation}
where $r = \sqrt{x_1^2 + x_2^2}$ and $\int\hpsngabs$ is defined as before.
\subsection{Nonlocal PDEs}
In many scientific and engineering problems, standard local models are not sufficient to accurately describe certain nonlocal phenomena, e.g., interactions at a distance. Hence nonlocal PDEs, which can express a more general description of the dynamical system, have been developed. Here we mention that peridynamics model for continuum mechanics \cite{silling2000reformulation,silling2005peridynamic} and  anomalous diffusion models \cite{du2012analysis,d2017nonlocal}.

In general, given the bounded, open domain $\Omega\subset \mathbb{R}^d$ and given a constant $\delta >0$, we define the \textit{interaction domain} corresponding to $\Omega$ as

\begin{equation}
	\Omega_{I_\delta}\coloneqq \{{y}\in \mathbb{R}^d\backslash\Omega \mbox{ such that } {y}\in B_{\delta}({x}) \mbox{ for some }{x}\in \Omega \},
\end{equation}
where $B_\delta({x})$ denotes the ball of radius $\delta$ centered at ${x}$.
For $\delta>0$, we consider the nonlocal problem \cite{d2020numerical} for a scalar-valued function $u({x})$ defined on $\Omega \cup \Omega_{I_\delta}$, given by
\begin{equation}
	\left\{
	\begin{array}{ll}
		-\mc{L}_{\delta}u(x)=f({x}),	&\forall {x}\in\Omega,\\
		\mc{V}u(x) = g({x}),&\forall {x}\in \Omega_{I_{\delta}}.
	\end{array}
	\right.
	\label{nonolocal_equation}
\end{equation}
Here $f({x})$ and $g({x})$ are given scalar-valued functions and
\begin{equation}
	\mc{L}_{\delta}u({x})\coloneqq 2\int _{\Omega\cup \Omega_{I_\delta}} (u({y})-u({x}))\gamma_\delta({x},{y})\mathrm{d}{y}\quad \mbox{for all } {x}\in \Omega,
	\label{nonlocal_operator}
\end{equation}
where $\gamma_\delta({x},{y})$ is a symmetric function, that is,
\begin{equation}
	\gamma _\delta({x}, {y})=\gamma_\delta ({y}, {x}),
\end{equation}
and for any ${x}$,
\begin{equation}
	\supp \big(\gamma_\delta({x},{y})\big)=B_\delta({x}).
\end{equation}

We assume $\gamma_\delta({x},{y})$ can be written in the form
\begin{equation}
	\gamma_\delta({x},{y})=\phi_\delta({x},{y})\theta_\delta({x},{y})\mc{X}_{B_\delta({x})}({y}),
\end{equation}
where $\theta_\delta({x},{y})$ and $\phi_\delta({x},{y})$ denote non-negative, symmetric, scalar-valued functions. We will refer to
$\gamma_\delta({x},{y})$ as the kernel, $\phi_\delta({x},{y})$ as the kernel function and $\theta_\delta({x},{y})$ as a constitutive function.
Because $\theta_\delta({x},{y})$ is a constitutive function which is not specific even within a single application, we focus on choices for the function $\phi_\delta({x},{y})$:
\begin{itemize}
	\item \textit{Translation-invariant, integral kernel functions}. $\phi_\delta({x},{y})=\phi_\delta({x}-{y})$ and satisfies for some positive constant $C>0$,
	$$C\leq \int _{\Omega\cup \Omega_{I_\delta}}\phi_\delta ({{y}}-{{x}})\mathrm{d}{{y}}<\infty \quad \mbox{for all }{x}\in \Omega\cup\Omega_{I_\delta}.$$
	\item \textit{``Critical'' kernel functions.}
	$$\phi_\delta ({x},{y})\propto \frac{1}{\Vert{y}-{x}\Vert^d}.$$
	\item \textit{``Peridynamic'' kernel functions}.
	$$\phi_\delta({x},{y})\propto \frac{1}{\Vert {y}-{x}\Vert}.$$
	\item \textit{Fractional kernel functions.} $$\phi_\delta({{x}}-{{y}}) \propto\frac{1}{\Vert{{y}}-{{x}}\Vert^{d+2s}},$$
	where $s\in(0,1).$
\end{itemize}
For more details one can refer to  \cite{d2020numerical} and references therein.

\subsection{Monte Carlo sampling for nonlocal operators}
Note that PINNs expresses the derivatives via automatic differentiation, which is not valid for IDEs and nonlocal operators. In this section, we propose a MC procedure to circumvent this situation. To this end, we first consider the IDE (\ref{integro_differential_equation}). For simplicity, we set $h_1(x)=0, h_2(x)=x$,
\begin{equation}
	u(x) = f(x) + \int_0^x K(x,s)u(s)ds.
\end{equation}
We mainly consider two cases:
\begin{itemize}
	\item Bounded kernel.
	If $K(x,s)$ is uniformly bounded for $x$ and $s$, then we adopt stochastic approximation via MC sampling as follows:
	\begin{equation}
		\int_0^xK(x,s)u(s)ds = x\cdot \mathbb{E}_{\xi\sim \mc{U}[0,1]} \left[K(x, x\xi) u(x\xi)\right],
	\end{equation}
	where $\xi$ is uniformly distributed on the interval $[0,1]$.
	\item Weakly singular kernel. Here, we set $K(x,s)=(x-s)^{-\alpha}, 0<\alpha<1$. In fact, it can be viewed as a generalized Abel's type equation, which occurs in many branches of scientific fields, such as microscopy, seismology, electron emission, plasma diagnostics.
	
	Notice that the kernel $|K(t,s)| = (t-s)^{-\alpha}$ is unbounded, yet integrable, hence one can view it as a scaled probability density function of Beta distribution. More precisely, we have
	\begin{equation}
		\int_0^x (x-s)^{-\alpha}u(s)ds=\int_0^x s^{-\alpha}u(x-s)ds=\frac{x^{1-\alpha}}{1-\alpha}\mathbb{E}_{\xi\sim Beta(1-\alpha,1)}\left[u(x-x\xi)\right].
	\end{equation}
	
	\item Hypersingular integral (\ref{hypersingular_def}). Using the definition (\ref{hypersingular_limit_def}), one can rewrite it as follows:
	\begin{equation}
		\begin{aligned}
			\int_a^b\hpsngAbs \frac{u(x)}{(x-s)^2}\mathrm{d}x=&\int_a^b\int_0^1(1-t)\frac{\partial ^2 u}{\partial x^2}\big((1-t)s+tx\big)\mathrm{d}t\,\mathrm{d}x\\
			&-u(s)\bigg(\frac{1}{b-s} + \frac{1}{s-a}\bigg) + \frac{\partial u}{\partial x}(s)\ln\frac{b-s}{s-a}.
		\end{aligned}
	\end{equation}
	Then we can approximate the above hypersingular integral via MC sampling
	\begin{equation}
		\begin{aligned}
				\int_a^b \hpsngAbs \frac{u(x)}{(x-s)^2}dx =& \mathbb{E}_{t\sim \mathcal{U}[0,1], x\sim \mathcal{U}[a,b]}\left[(1-t)\frac{\partial ^2 u}{\partial x^2}\big(tx + (1-t)s\big)\right]\\ &-u(s)\bigg(\frac{1}{b-s} + \frac{1}{s-a}\bigg) + \frac{\partial u}{\partial x}(s)\ln\frac{b-s}{s-a}.
		\end{aligned}
	\end{equation}
	For the hypersingular integral (\ref{two_dim_hyper_singular}), one can rewrite it as follows:
	\begin{equation}
		\begin{aligned}
		&{\hiint}_\Omega \frac{u(x_1, x_2)}{r^3}dx_1dx_2\\ &=\int_0^{2\pi}\left[-\frac{1}{R(\nu)}\tilde{u}(0,\nu) + \frac{\partial \tilde{u}}{\partial r}(0,\nu)\ln(R(\nu)) + \int_0^{R(\nu)}\int_0^1(1-t)\frac{\partial ^2 \tilde{u}}{\partial r^2}(tr,\nu)dtdr\right]d\nu,
		\end{aligned}
	\end{equation}
	where $\tilde{u}(r,\nu) = u(r\cos \nu, r\sin \nu).$ Then we approximate it via  MC sampling
	\begin{equation}
		\begin{aligned}
					&{\hiint}_\Omega \frac{u(x_1, x_2)}{r^3}dx_1dx_2 \\&= \mathbb{E}_{\nu\sim \mc{U}[0,2\pi]}\left[-\frac{1}{R(\nu)}\tilde{u}(0,\nu) + \frac{\partial \tilde{u}}{\partial r}(0,\nu)\ln(R(\nu)) + \mathbb{E}_{r\sim \mc{U}[0,R(\nu)],t\sim \mc{U}[0,1]} (1-t) \frac{\partial ^2 \tilde{u}}{\partial r^2}(tr,\nu)
			\right].
		\end{aligned}
	\end{equation}
	
	For the nonlocal model (\ref{nonolocal_equation}) with the integral operator (\ref{nonlocal_operator}), our goal is to approximate
	\begin{equation}
		\mathcal{L}_{\delta}u(x)=2\int_{\Omega\cup \Omega_\delta} (u(y)-u(x))\gamma_\delta(x,y)dy.
	\end{equation}
	We assume that \begin{equation}
		\gamma_\delta(x,y) = \frac{1}{\Vert y-x\Vert ^{\alpha}}\cdot \mathbbm{1}_{\Vert y-x\Vert < \delta},\quad \alpha \in (0,d+2).
	\end{equation}
	The corresponding stochastic approximation is
	\begin{equation}
		\begin{aligned}
			\mathcal{L}_\delta u(x) & = 2\int_{y\in B_\delta(x)}\frac{u(y)-u(x)}{\Vert y-x\Vert_2^{\alpha}} dy = 2\int _{\Vert y\Vert_2 <\delta}\frac{u(x+y)-u(x)}{\Vert y\Vert _2^\alpha}dy\\
			&=\int_{\Vert y\Vert_2<\delta}\frac{u(x+y)-2u(x)+u(x-y)}{\Vert y\Vert _2^\alpha}dy\\
			& = \int_{S^{d-1}}\int_0^\delta \frac{u(x+r\xi)-2u(x)+u(x-r\xi)}{r^\alpha} \cdot r^{d-1}drd\xi\\
			&=\left\{
			\begin{array}{lr}
				\delta ^{d-\alpha}|S^{d-1}|\cdot 	\mathbb{E}_{\xi, r\sim U[0,1]}\left[\frac{u(x+\delta r\xi)-2u(x)+u(x-\delta r\xi)}{r^{\alpha+1-d}}\right],&\alpha\in (0,d],\\
				\frac{\delta^{d+2-\alpha}}{d+2-\alpha}|S^{d-1}|\cdot 	\mathbb{E}_{\xi, r\sim Beta(d+2-\alpha,1)}\left[\frac{u(x+\delta r\xi)-2u(x)+u(x-\delta r\xi)}{\delta ^2 r^{2}}\right],&\alpha\in (d, d+2),
			\end{array}
			\right.
		\end{aligned}
	\end{equation}
	where $\xi$ is uniformly distributed on $S^{d-1}$, $|S^{d-1}|$ denotes the surface area of $S^{d-1}$.
	
	Note that as $r\to 0$ we have
	\begin{equation}
		\lim_{r\to 0^+}\frac{u(x+\delta r\xi)-2u(x)+u(x-\delta r\xi)}{\delta ^2 r^{2}} = \frac{\partial _r^2u(x)}{\partial r^2}\bigg|_{r=0},
	\end{equation}
	which may suffer from rounding errors for extremely small $\delta r$. Thus we truncate $r$ with $r_\epsilon = \max\{r,\epsilon\}$, and replace  variable $r$ in the integrand with $r_\epsilon$, where $\epsilon>0$ is a small positive number.
	
	%\begin{remark}
	%	In fact, for any singular but integrable kernel function, the above procedure is also true. Specifically, if the kernel of interest has more than one singularity, we can divide the domain (similarly with (\ref{laplacian_two_parts})) into many subdomains such that each subdomain has only one singularity. If the behavior near singularity(here we take $0$) of the kernel function does not have the expression like $K(x)=x^{-\alpha}$, one can at least modify the kernel as follows:
	%	\begin{equation}
		%		\widetilde{K}(x) = x^{-\alpha_0}\cdot \left(K(x)x^{\alpha_0}\right),
		%	\end{equation}
	%	where $\alpha_0\in(0,1)$ and $\left|K(x)x^{\alpha_0}\right|<\infty$.
	%	
	%\end{remark}
\end{itemize}
\subsection{MC-Nonlocal-PINNs}
\begin{figure}[H]
	\centering
	\includegraphics[scale=0.55]{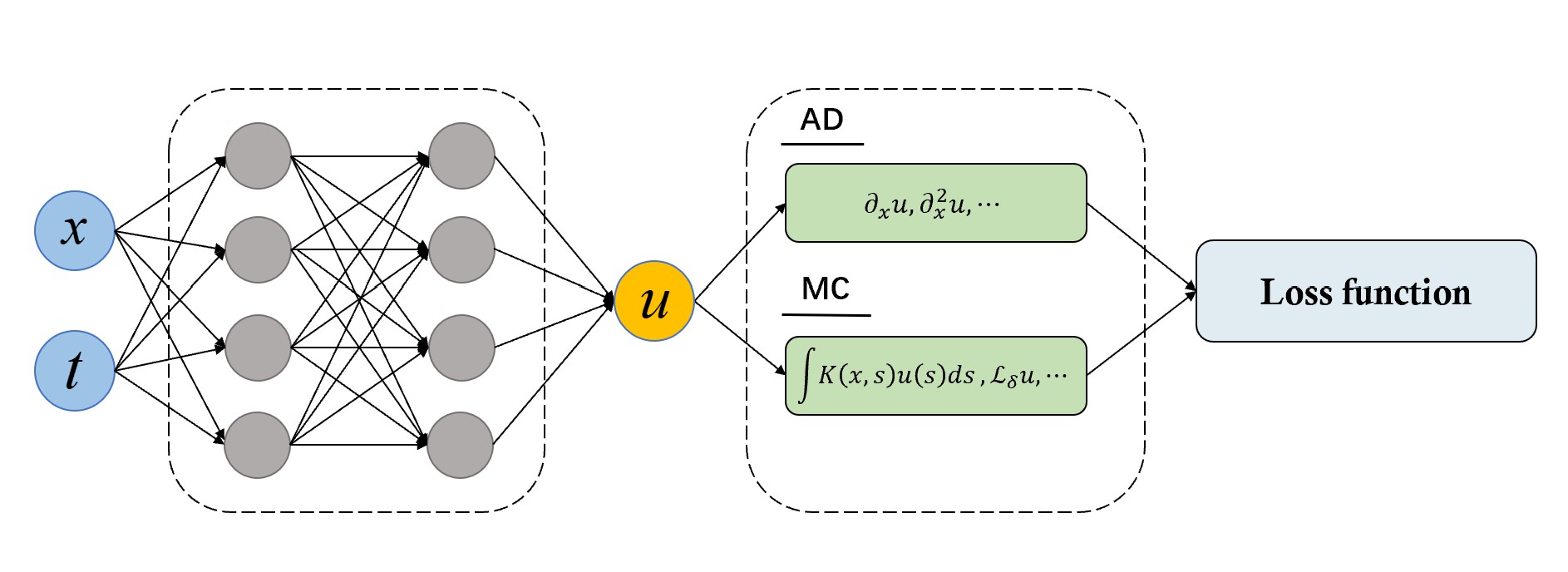}
	\caption{The architecture of MC-Nonlocal-PINNs.}
	\label{fig:integralarchitecture}
\end{figure}

Based on the above stochastic approximation, we are ready to present our algorithm.
For IDEs, because the equation in (\ref{integro_differential_equation}) combines the differential operator and the integral operator, it is necessary to define boundary conditions
\begin{equation}
	u(x) = g(x),\quad x\in \Gamma.
\end{equation}
For the nonlocal model, we must address the corresponding non-zero volume boundary $u(x)(x\in \Omega_{I_\delta})$. In neural network framework, we address it by adding additional soft penalties in the final loss function. The schematic of MC-Nonlocal-PINNs is shown in Figure \ref{fig:integralarchitecture}. Specifically, we have access to some collocation points $\{x_r^i\}_{i=1}^{N_r} \subset \Omega, \{x_b^i\}_{i=1}^{N_b}\subset \partial \Omega$ (or $\Omega_{I_\delta}$). Following the stochastic approximation and automatic differentiation, we can calculate the loss functions $\mathcal{L}_r(\theta),\mathcal{L}_b(\theta)$. The weights for each component are given by adaptive strategy:
\begin{equation}
	[w_r,w_b] = \frac{{[\mathcal{L}_r(\theta),\mathcal{L}_b(\theta)]}}{\min\{\mathcal{L}_r(\theta),\mathcal{L}_b(\theta)\}}.
	\label{weight}
\end{equation}
We summarize our method in Algorithm \ref{algorithm}.
\begin{algorithm}[H]
	\begin{itemize}
		\item\textbf{1.} Specify the training set
		\begin{equation*}
			\mathcal{D}  = \left\{\{x_r^i\}_{i=1}^{N_r}, \{x_b^i\}_{i=1}^{N_b}\right\}.
		\end{equation*}
		\item\textbf{2.} Sample $N$ snapshots from the above training data
		\item\textbf{3.} Calculate the loss $\mathcal{L}(\theta)=w_r\cdot \mathcal{L}_r(\theta) +w_b \cdot \mathcal{L}_b(\theta)$ for  via (\ref{loss}) and (\ref{weight})
		\item\textbf{4.} Let $W\leftarrow W-\eta\frac{\partial \mathcal{L}}{\partial W}$ to update all the involved parameters $W$, where $\eta$ is the learning rate
		\item\textbf{5.} Repeat \textbf{Steps 2-4} until convergence
	\end{itemize}
	\caption{MC-Nonlocal-PINNs}
	\label{algorithm}
\end{algorithm}

\section{Numerical experiments}
In this section, we present a series of comprehensive numerical tests to demonstrate the effectiveness of proposed algorithm. We investigate the performance of MC-Nonlocal-PINNs for solving Volterra-type equations and hypersingular integral equations, then we illustrate the efficiency of the MC-Nonlocal-PINNs method to solve general nonlocal equations. To quantitatively evaluate the accuracy of numerical solution, we shall consider $L^2$ relative error of the predicted solution:
$$\mbox{Relative $L^2$ error} = \frac{\Vert u_{NN}(x) - u(x)\Vert_2}{\Vert u(x)\Vert_2},$$
where $u$ and $u_{NN}$ are fabricated and surrogate solutions, respectively.

Throughout all experiments, the DNNs model contains four hidden layers with 64 neurons per hidden layer. We shall employ hyperbolic tangent activation functions (Tanh) and initialize all trainable parameters using Glorot initialization, unless stated otherwise. All networks are trained using the Adam optimizer with default settings and the L-BFGS optimizer. We adopt exponential learning rate decay with a decay-rate of 0.9 every 1000 training iterations.
\subsection{Volterra integral equation}
\subsubsection{1D bounded kernel problem}
\begin{figure}[!ht] %并排插入多个图片
	\begin{minipage}[t]{0.32\textwidth}
		\centering
		\includegraphics[scale=0.3]{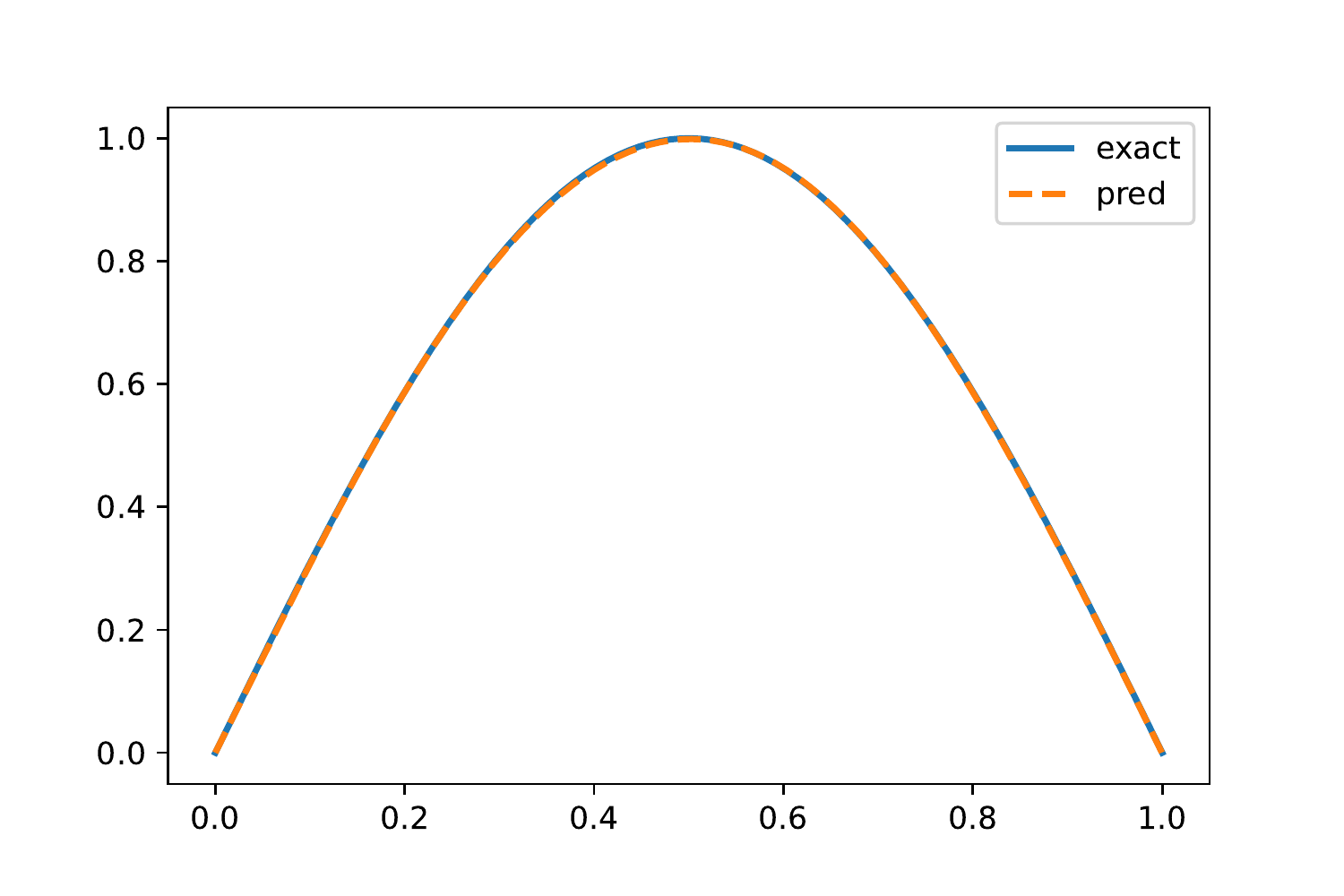}
	\end{minipage}
	\begin{minipage}[t]{0.32\textwidth}
		\centering
		\includegraphics[scale=0.3]{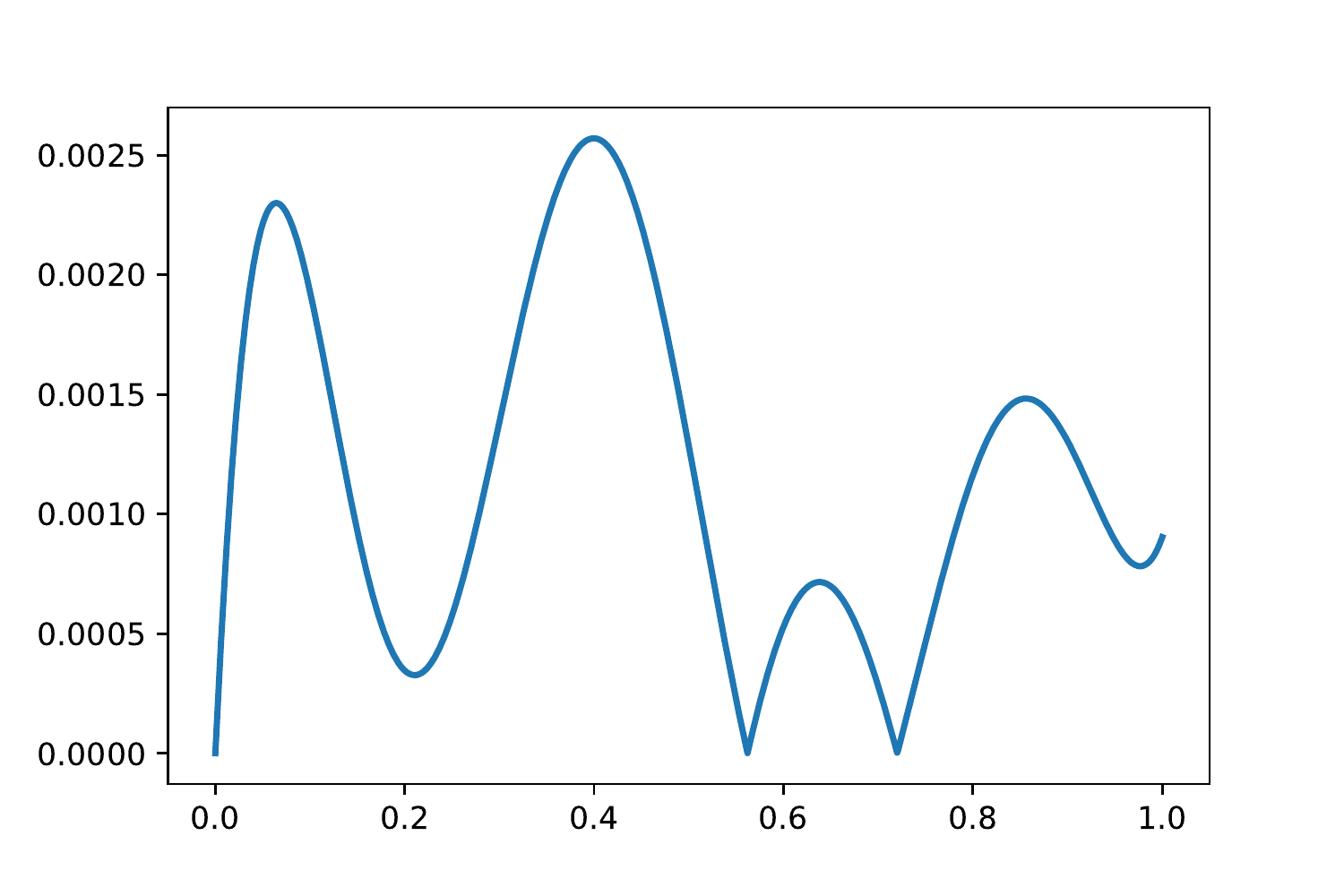}
	\end{minipage}
	\begin{minipage}[t]{0.32\textwidth}
		\centering
		\includegraphics[scale=0.3]{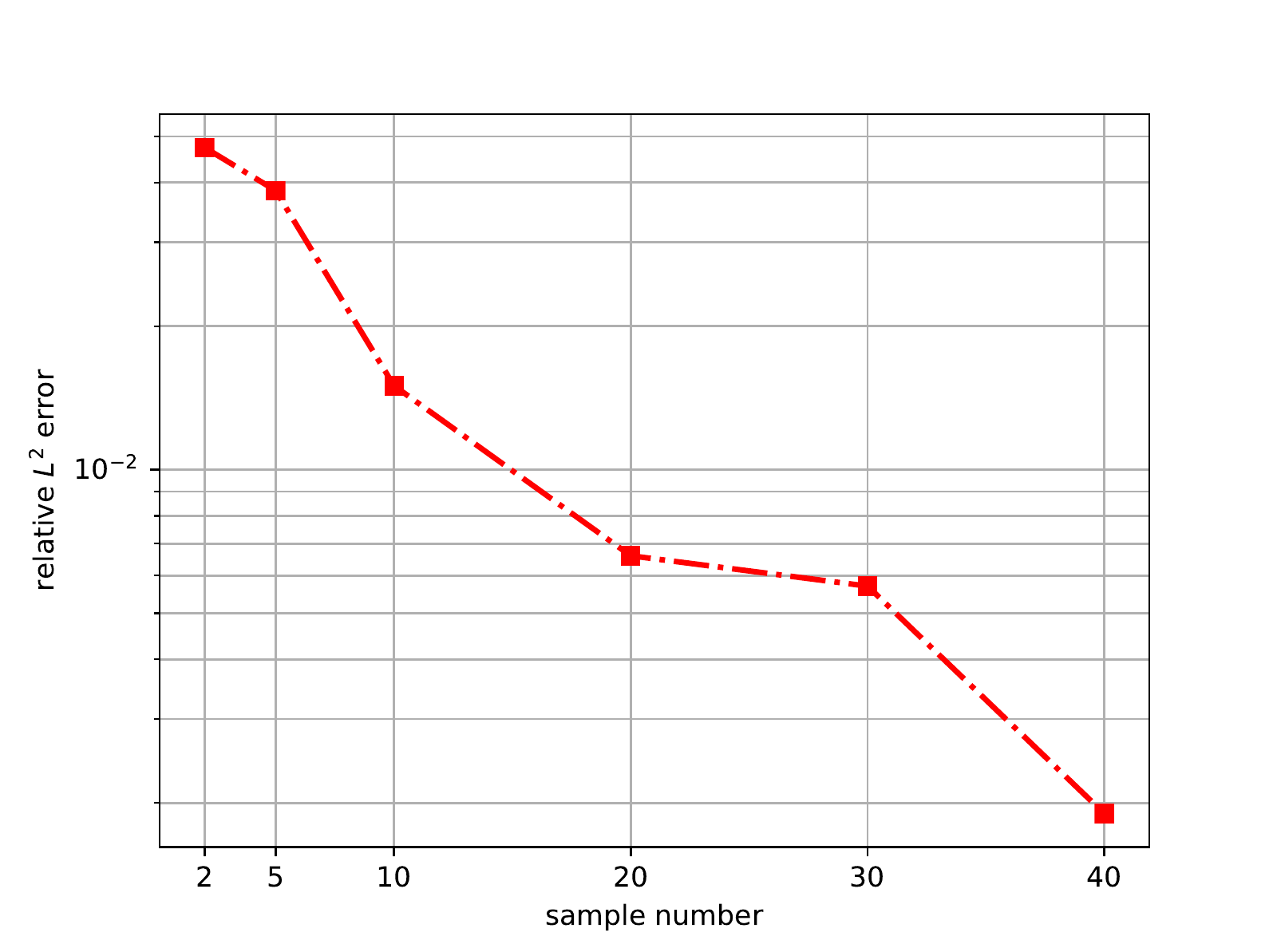}
		\label{fig:volterra_1D_bounded_kernel_convergence}
	\end{minipage}
	\caption{Volterra equation (1D bounded kernel). Left : exact and predicted solutions. Middle : absolute error. Right : the relative $L^2$ error for different $N_s$.}
	\label{fig:volterra_1D_bounded_kernel_solution}
\end{figure}

Consider the following example:
\begin{equation}\label{1D_volterra_eqn}
	u(x) = f(x) + \int_{0}^x K(x,s)u(s)ds, \quad 0\leq x\leq 1,
\end{equation}
where the exact solution and the corresponding terms are given as follows:
\begin{equation*}
	u(x) = \sin(\pi x), \quad f(x) = \left(1-\frac{1}{2\pi}\right)\sin(\pi x)-\cos(\pi x) / {2\pi}, \quad K(x,s)=-\sin(\pi (x-s)).
\end{equation*}
Note that kernel $K(x,s)$ is bounded in $[0,1]\times [0,1]$, one can rewrite the equation as follows
\begin{equation*}
	u(x) = f(x) + \mathbb{E}_{s\sim U[0,1]}\left[x\cdot K(x,xs)\cdot u(xs)\right].
\end{equation*}
As is described in Section \ref{section:methodology}, we approximate the expectation using Monte Carlo method:
$$\mathbb{E}_{s\sim U[0,1]}\left[x \cdot K(x,xs)\cdot u_{NN}(xs;\theta)\right]\approx \frac{1}{N_s}\sum_{i=1}^{N_s}x \cdot K(x, xs_i)\cdot u_{NN}(xs_i;\theta),$$
where $s_i \sim U[0,1]$ and $N_s$ is the number of discrete integration points. We use 128 uniformly
distributed training points in the space domain for each batch and train the MC-Nonlocal-PINNs using Adam optimizer with an initial learning rate of 0.001 to $1000$ iterations, then we continue to train the model using L-BFGS with adaptive learning rate to $1000$ iterations. Figure \ref{fig:volterra_1D_bounded_kernel_solution} shows the comparison between the predicted and the exact solution  for $Ns=40$. We observe that the predictions achieve an excellent agreement with the corresponding ground truth. Furthermore, we investigate the performance under the cases of different sample numbers $N_s$, the final result is reported on the right of Figure \ref{fig:volterra_1D_bounded_kernel_solution}. When increasing the sample number $N_s$, the relative $L^2$ error of $u_{NN}$ decreases from 4.74\% to 0.19\%.

\subsubsection{1D weakly singular kernel problem}
\begin{figure}[H] %并排插入多个图片
	\begin{minipage}[t]{0.32\textwidth}
		\centering
		\includegraphics[scale=0.3]{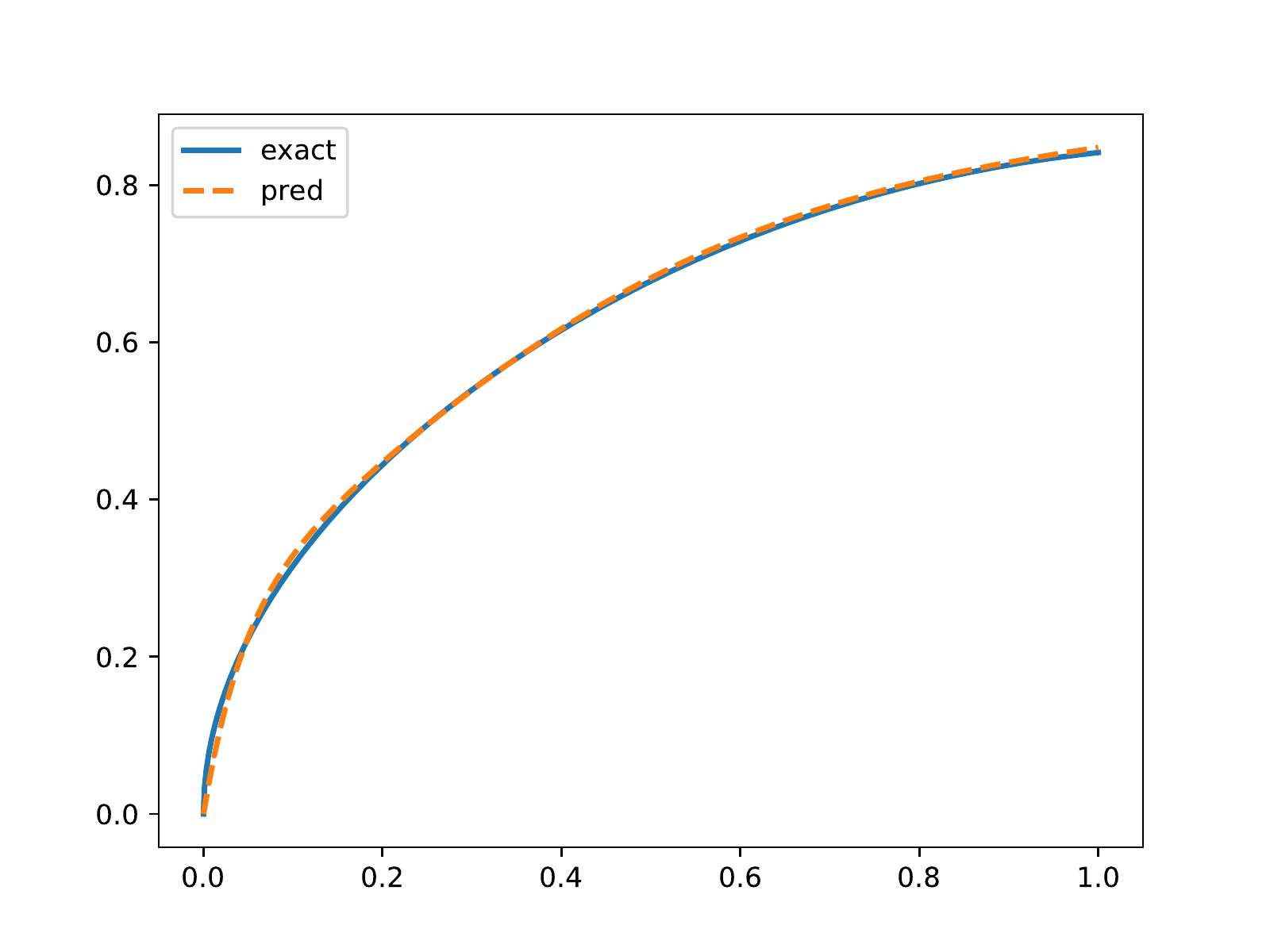}
	\end{minipage}
	\begin{minipage}[t]{0.32\textwidth}
		\centering
		\includegraphics[scale=0.3]{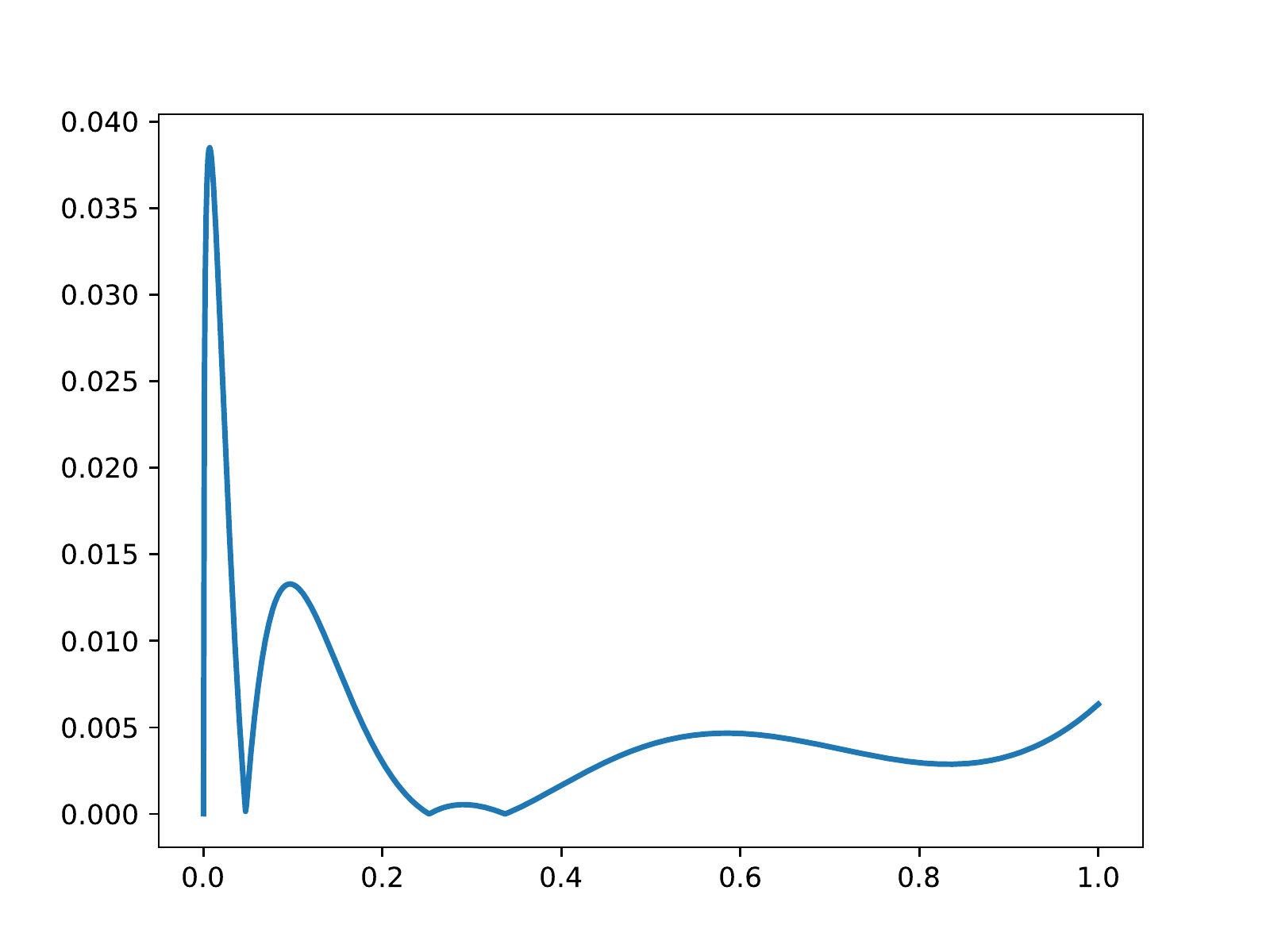}
	\end{minipage}
	\begin{minipage}[t]{0.32\textwidth}
		\centering
		\includegraphics[scale=0.3]{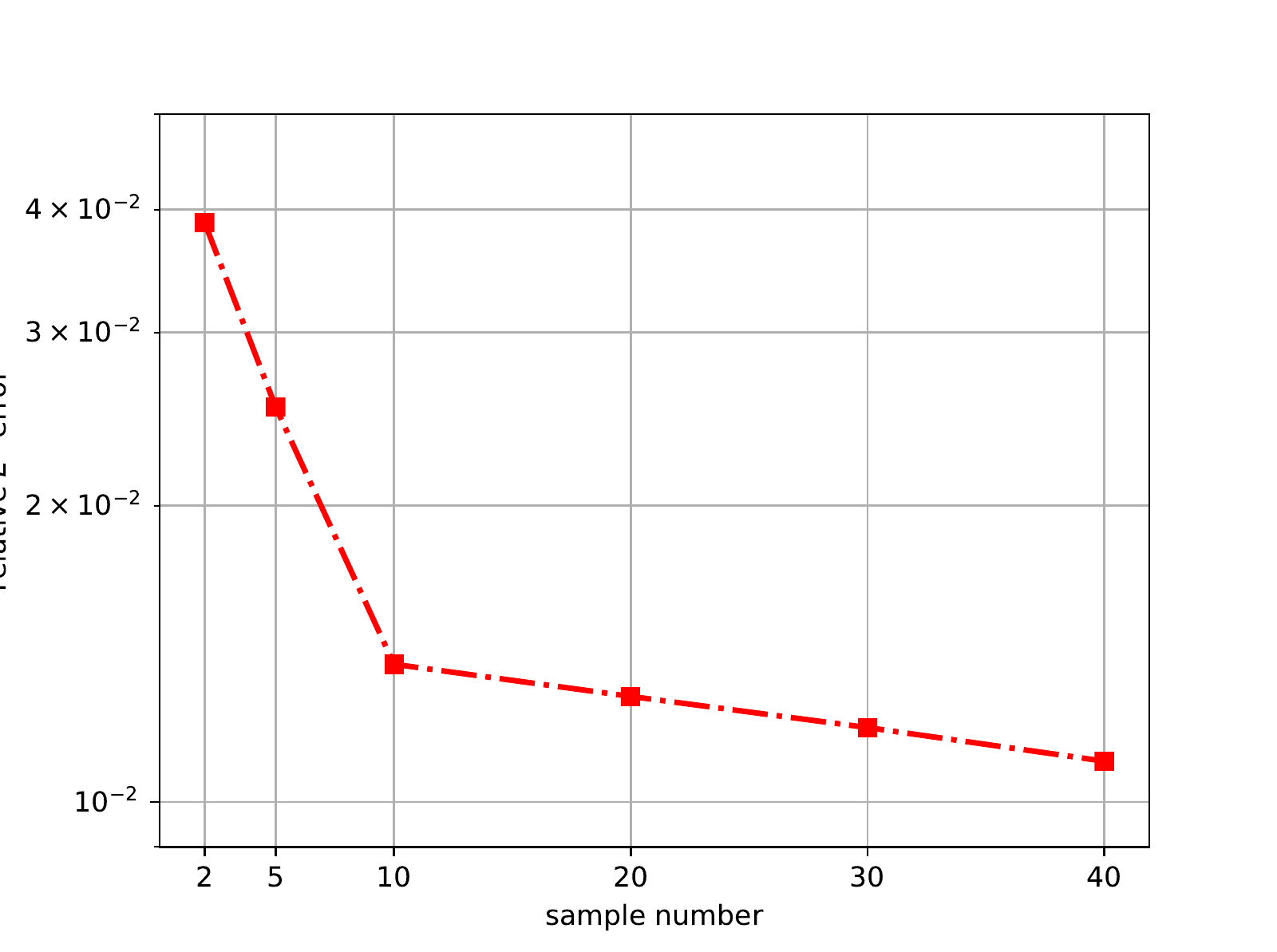}
	\end{minipage}
	\caption{Volterra equation (1D singular kernel). Left : exact and predicted solutions. Middle : absolute error. Right : the relative $L^2$ error for different $N_s$.}
	\label{fig:volterra_1D_singular_kernel_solution}
\end{figure}

In this example, we still consider the above equation (\ref{1D_volterra_eqn}) with a singular kernel. Specifically, we take
\begin{equation*}
	u(x)=\frac{\sin (x)}{\sqrt{x}},\quad f(x) = \frac{\sin(x)}{\sqrt{x}}+\pi \sin \left(\frac{x}{2}\right)J_0\left(\frac{x}{2}\right),\quad K(x,s)=-(x-s)^{-\alpha},
\end{equation*}
where $\alpha=1/2$ and $J_0(z)$ is the Bessel function of the first kind defined by
\begin{equation}
	J_0(z)=\sum_{k=0}^{\infty}\frac{(-z^2)^k}{(k!)^2 4^k}.
\end{equation}
Similarly, we can rewrite the equation as
\begin{equation*}
	u(x) = f(x) + 2\sqrt{x}\cdot \mathbb{E}_{s\sim Beta(0.5,1)}\left[u(x-xs)\right],
\end{equation*}
and approximate the expectation using MC sampling
\begin{equation*}
	u_{NN}(x;\theta) \approx f(x) + \frac{2\sqrt{x}}{N_s}\cdot \sum_{i=1}^{N_s} u_{NN}(x -xs_i;\theta),
\end{equation*}
where $s_i \sim Beta(0.5,1)$ and $N_s$ is the number of discrete integration points. We set batch size to 128 and $N_s$ to 100, and train the MC-Nonlocal-PINNs using the Adam optimizer with initial learning rate 0.001 to $1000$ iterations, then we continue to train the model using L-BFGS with adaptive learning rate to $2000$ iterations. The results for the exact and the predicted solutions are presented in Figure \ref{fig:volterra_1D_singular_kernel_solution}. As is shown, a good agreement can be achieved between the ground truth and predicted solution. And relative $L^2$ errors for different sample number $N_s$ are shown on the right of Figure \ref{fig:volterra_1D_singular_kernel_solution}.

\subsubsection{1D Fredholm problem}
\begin{figure}[H] %并排插入多个图片
	\begin{minipage}[t]{0.32\textwidth}
		\centering
		\includegraphics[scale=0.3]{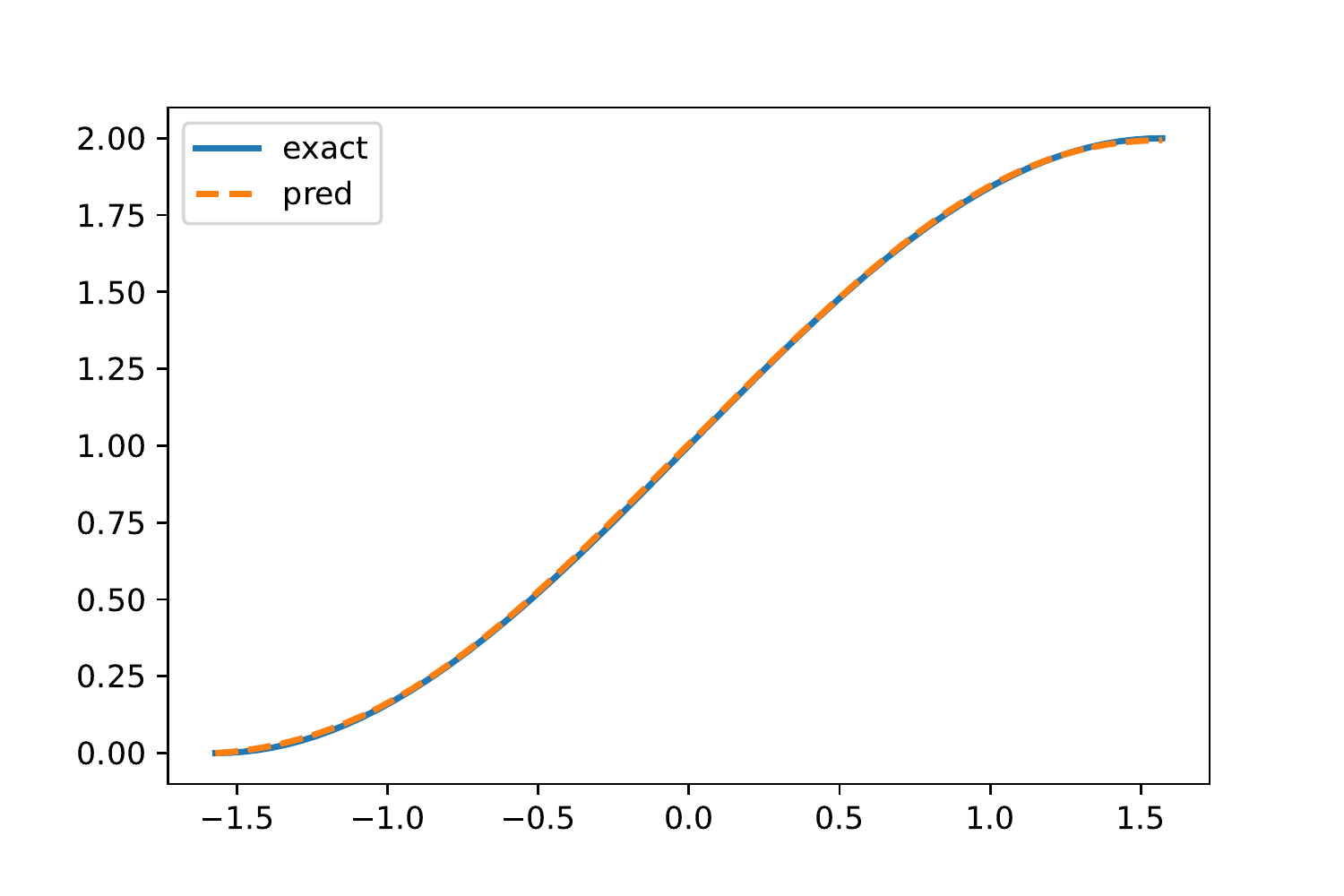}
	\end{minipage}
	\begin{minipage}[t]{0.32\textwidth}
		\centering
		\includegraphics[scale=0.3]{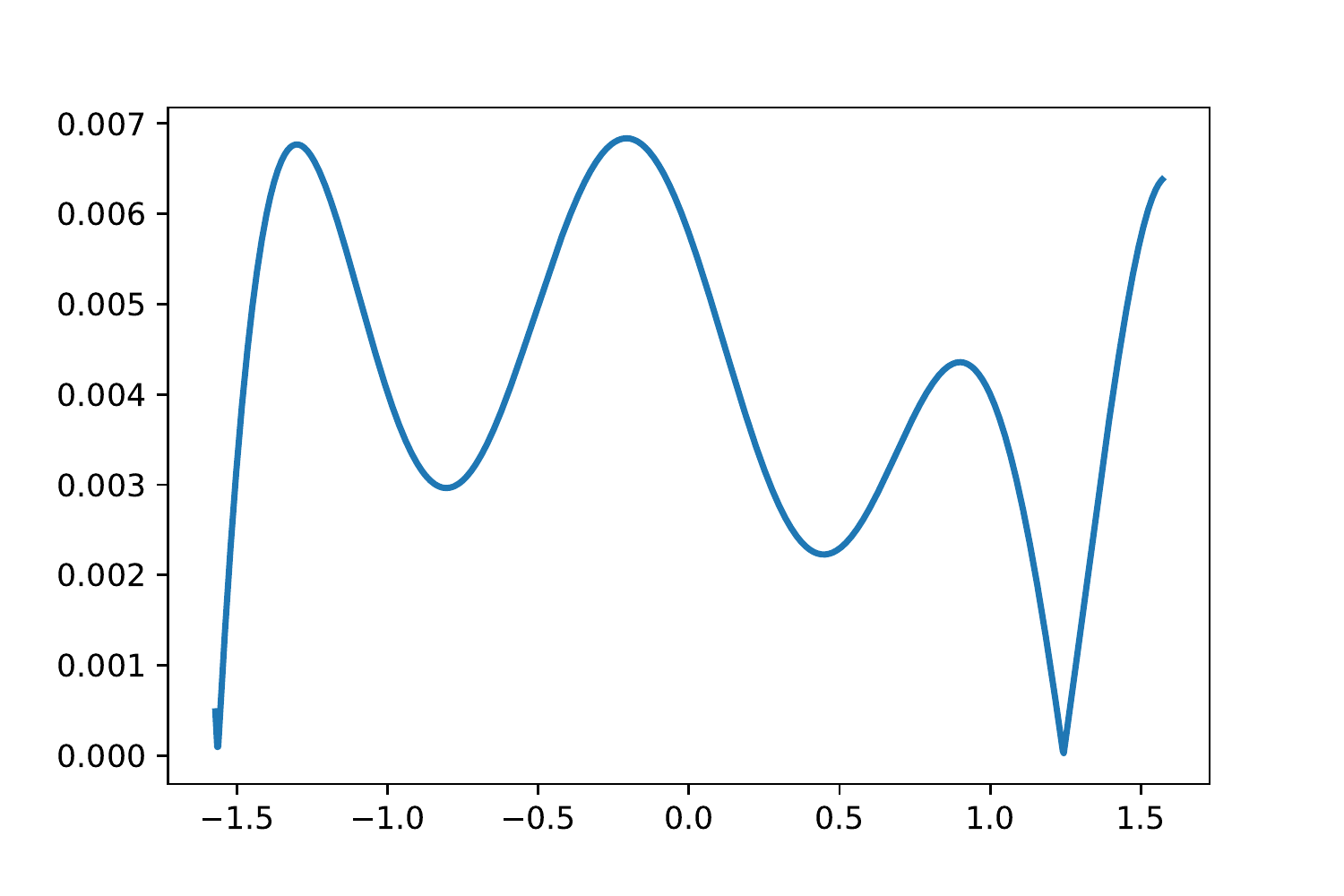}
	\end{minipage}
	\begin{minipage}[t]{0.32\textwidth}
		\centering
		\includegraphics[scale=0.3]{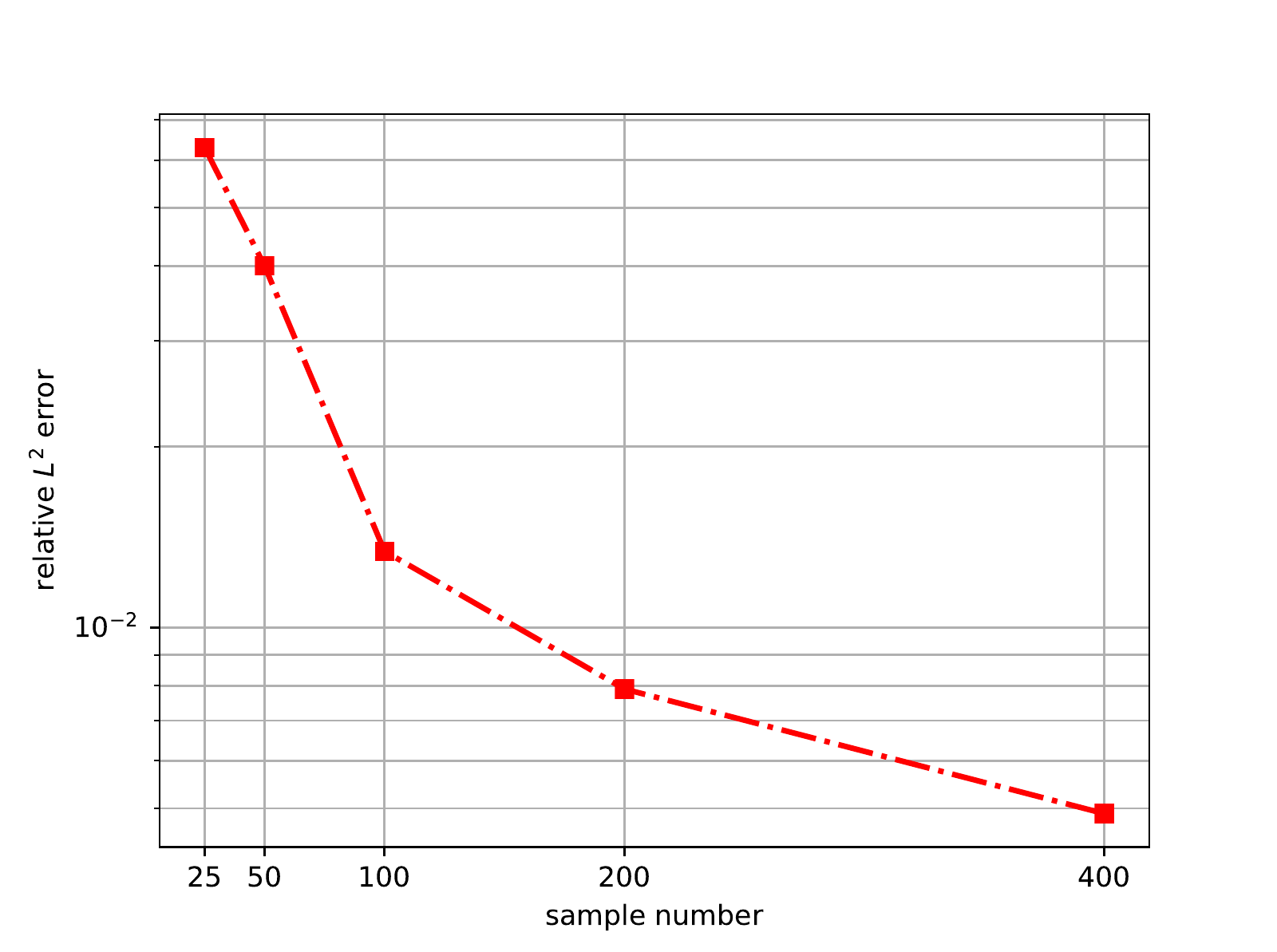}
	\end{minipage}
	\caption{Fredholm equation. Left : exact and predicted solutions. Middle : absolute error. Right : the relative $L^2$ error for different $N_s$.}
	\label{fig:Fredholm_solution}
\end{figure}

In this section, we consider a nonlinear 1D Fredholm IDE:
\begin{equation*}
	\frac{du}{dx}=\cos (x) - x + \frac{1}{4}\int_{-\pi/2}^{\pi/2}xtu^2(t)dt,\quad u(-\frac{\pi}{2})=0.
\end{equation*}
And the corresponding exact solution is chosen as $u(x)=1+\sin(x)$. Since the limits of integration are constants, we can approximate the expectation using uniform sampling method
\begin{equation*}
	\frac{\partial u_{NN}(x;\theta)}{\partial x} \approx \cos(x) -x + \frac{\pi}{4N_s} \sum_{i=1}^{N_s} xs_i \cdot u_{NN}(s_i;\theta)^2,
\end{equation*}
where $s_i\sim U[-\pi/2,\pi/2]$. We use 128 uniformly
distributed training points in the space domain for each batch. We take $N_s=400$ and train the MC-Nonlocal-PINNs using the Adam optimizer with initial learning rate 0.001 to $5000$ iterations, then we continue to train the model using L-BFGS with adaptive learning rate to $2000$ iterations. And the final results are shown in Figure \ref{fig:Fredholm_solution}. We can observe that as sample number ($N_s$) increases, the relative $L^2$ error decreases from 6.3\% to 0.49\%.

\subsubsection{High dimensional bounded kernel problem}
Consider the following 10D Volterra IDE \cite{yuan2022pinn}:
\begin{equation*}
	\left\{
	\begin{aligned}
		&\frac{\partial u(t,x_1,\cdots,x_9)}{\partial t}+\sum_{i=1}^{9}\frac{\partial u(t,x_1,\cdots,x_9)}{\partial x_i}=f(t,x_1,\cdots,x_9),\\
		&f(t,x_1,\cdots,x_9) = u(t,x_1,\cdots,x_9) + g(t,x_1,\cdots,x_9)+\int_0^{x_9}\cdots\int_0^{x_1}\int_0^{t}s_0\cdot u(s_0,s_1,\cdots,s_9)ds_0ds_1\cdots ds_{9},
	\end{aligned}
	\right.
\end{equation*}
where $0\leq t,x_1,\cdots,x_9\leq 1.$ The exact solution is
\begin{equation*}
	u(t,x_1,\cdots,x_9)=t\cdot (x_1+x_2+x_3)\cdot \sin (x_4+x_5+x_6)\cdot \cos(x_7+x_8+x_9).
\end{equation*}

\begin{figure}[H]
	\centering
	\includegraphics[scale=0.4]{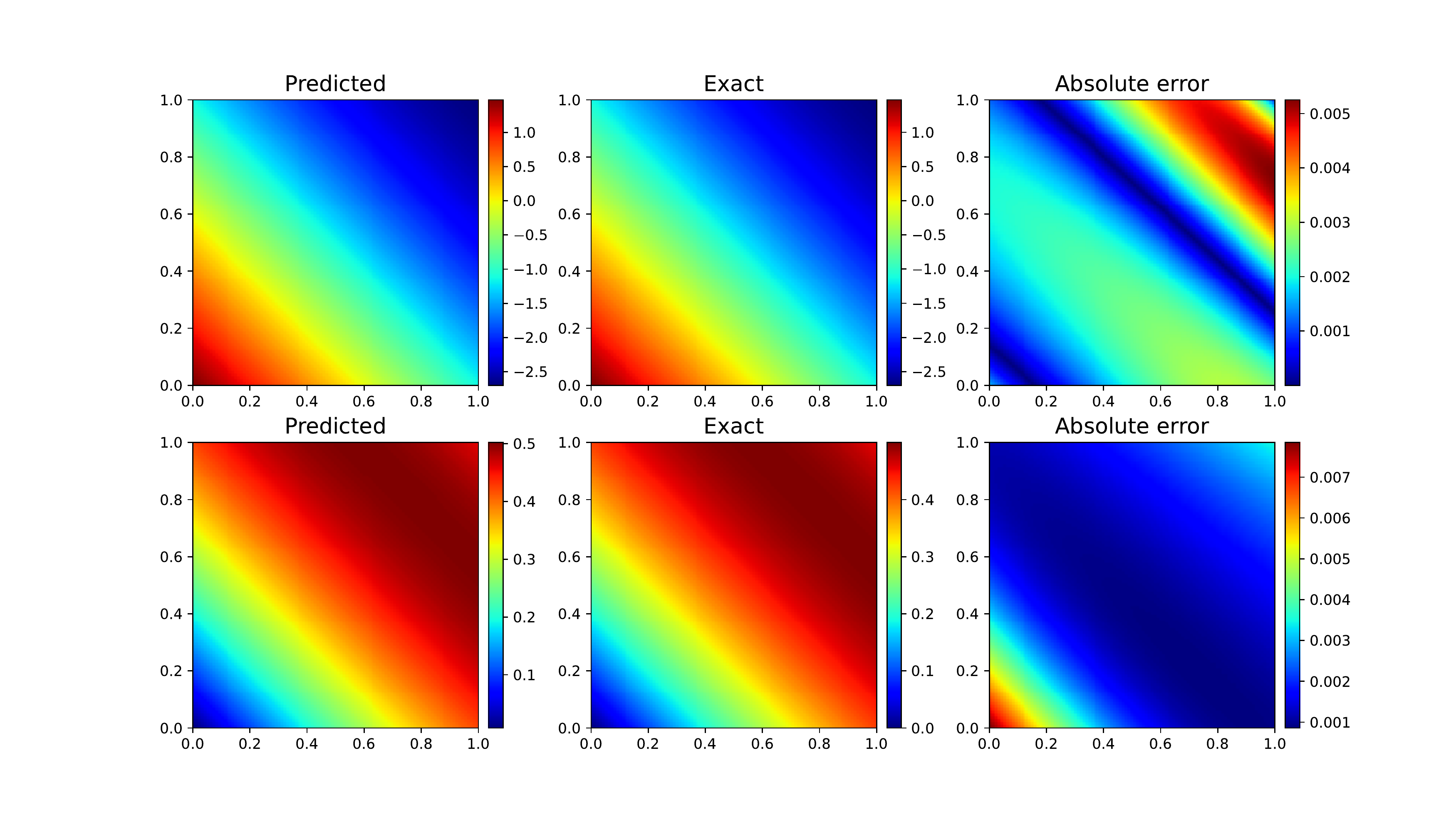}
	\caption{10D Volterra equation. From left to right: predicted solutions, exact solutions and corresponding absolute errors.}
	\label{fig:volterra_high_dim_bounded_kernel_result}
\end{figure}

Our goal is to approximate the 10-dimensional integral terms. After some simple calculations, we have
$$\int_0^{x_9}\cdots\int_0^{x_1}\int_0^{t}s_0\cdot u(s_0,s_1,\cdots,s_9)ds_0ds_1\cdots ds_{9}\approx x_1x_2\cdots x_9\cdot t^2\cdot\frac{1}{N_s} \sum_{i=1}^{N_s} u(ts_0^i, x_1s_1^i,\cdots,x_9s_9^i),$$
where $(s_0^i,s_1^i,\cdots,s_9^i)\sim U[0,1]^{10}.$ The training set consists of two parts: 10000 collocation points randomly sampled in the equation domain and 1000 boundary points randomly sampled on each boundary. We take 10 Monte Carlo points to approximate integral terms, i.e., $N_s=10$, and train the MC-Nonlocal PINNs using LBFGS optimizer 40000 iterations. To elucidate the solution of our method, we select two different planes $[1,1,1,1,0,1,1,x_8,x_9]$ and $[1,1,0,0,x_4,x_5,0,0,0,0]$. The corresponding exact solutions are
$u(0.5,1,1,1,0,1,1,x_8,x_9)=3\cdot \sin(2)\cdot \cos(x_8+x_9+1)$ and $u(0.5,1,0,0,x_4,x_5,0,0,0,0) = \sin(x_4+x_5)/2$. The final results are shown in Figure \ref{fig:volterra_high_dim_bounded_kernel_result} and final relative $L^2$ error is $0.134\%$. It is noticed that our method achieves better accuracy than the A-PINN approach (0.519\%) without introducing 10 auxiliary variables.

\subsubsection{High dimensional singular kernel problem}
%Consider the following PDE:
%\begin{equation*}
%	\begin{aligned}
	%			-\Delta u(x) + \int_{\Vert s\Vert _2\leq \Vert x\Vert _2} \Vert x-s\Vert _2 ^{-\alpha}u(s)ds = f(x), \quad x\in \Omega,\\
	%			u(x)=0, \quad x\in \partial \Omega.
	%	\end{aligned}
%\end{equation*}
%Here $\alpha \in (0, d), \Omega=\left\{x\mid \Vert x\Vert _2\leq 1\right\}\subset \mathbb{R}^d$,  $d$ is spatial dimension. The exact solution $$u(x)=\big(1-\Vert x\Vert_2^2\big)^{1+\alpha}.$$
%We have
%\begin{equation*}
%	\begin{aligned}
	%		\int_{\Vert s\Vert _2\leq \Vert x\Vert _2} \Vert x-s\Vert _2 ^{-\alpha}u(s)ds &= \int_{S^{d-1}}\int_0^{\Vert x\Vert_2} r^{d-1-\alpha}u(x-r\theta)drd\theta \\
	%		&=\left\{
	%		\begin{array}{cc}
		%			\Vert x\Vert_2 ^{d-\alpha}\cdot\mathbb{E}_{\xi,r\sim U[0,1]} \left[r^{d-1-\alpha}u(x-xr\theta)\right],& \alpha \leq d-1,\\
		%			\frac{\Vert x\Vert_2^{d-\alpha}}{d-\alpha}\cdot \mathbb{E}_{\xi, r\sim Beta(d-\alpha, 1)}\left[u(x-xr\theta)\right],  & d-1<\alpha<d,
		%		\end{array}
	%		\right.		
	%	\end{aligned}
%\end{equation*}
%where $\xi$ is uniformly distributed on the sphere $S^{d-1}$, $|S^{d−1}|$ denotes the surface area of $S^{d-1}$.
Consider the following PDE:
\begin{equation*}
	\begin{aligned}
		&\frac{\partial u(t,x_1,\cdots,x_d)}{\partial t} + \sum_{i=1}^{d}\frac{\partial u(t,x_1,\cdots,x_d)}{\partial x_i} = f(t, x_1,\cdots,x_d),\\
		&f(t,x_1,\cdots,x_d) = u(t, x_1,\cdots,x_d) + g(t,x_1,\cdots,x_d) \\
		&\qquad\qquad \qquad \quad + \int_{0}^{x_d}\cdots \int_0^{x_1}\int_0^t(t-s_0)^{-\alpha}(x_1-s_1)^{-\alpha}\cdots (x_d-s_d)^{-\alpha} u(s_0,s_1,\cdots,s_d) dt ds_1\cdots ds_d.
	\end{aligned}
\end{equation*}
Here $\alpha=1/2, \Omega=[0,1]^{d+1}$,  $d$ is spatial dimension.  The exact solution is $$u(x)=\big(1-\Vert x\Vert_2^2\big)e^{-t}.$$

\begin{figure}[!ht]
	\centering
	\includegraphics[scale=0.38]{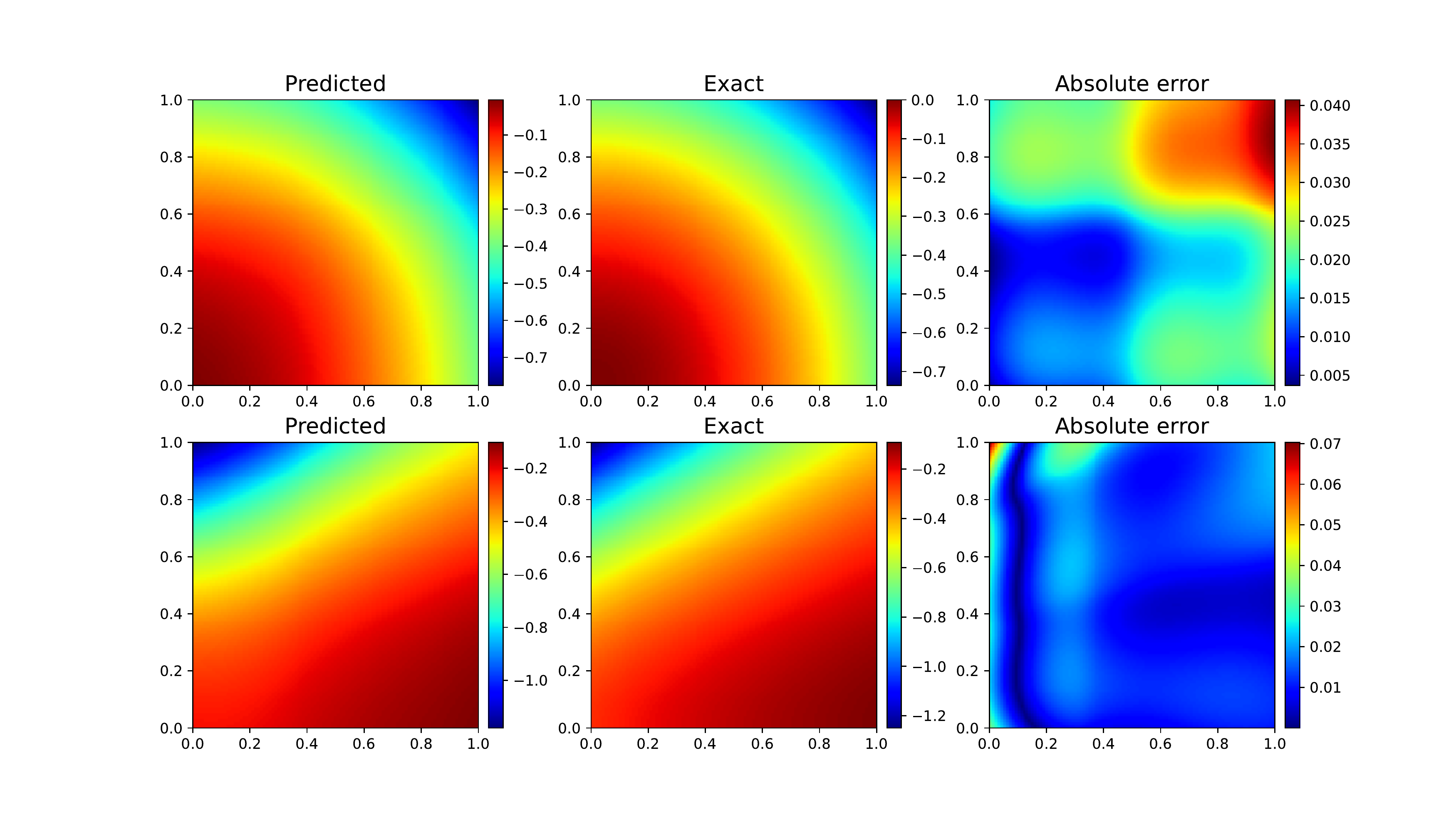}
	\caption{4D Volterra equation. From left to right: predicted solutions, exact solutions and the corresponding absolute errors.}
	\label{fig:volterra_4D_singular_kernel_result}
\end{figure}

\begin{figure}[!ht]
	\centering
	\includegraphics[scale=0.38]{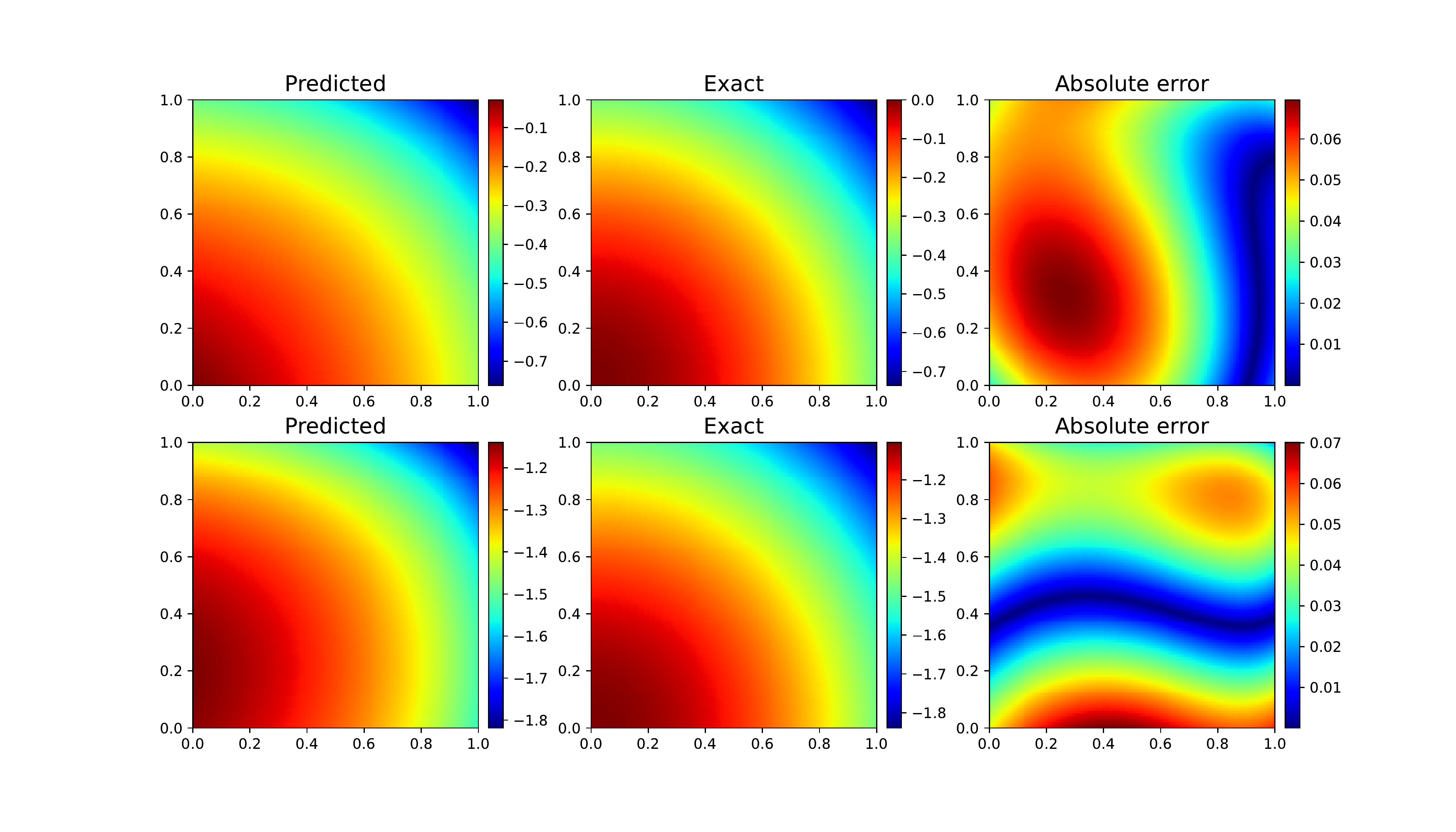}
	\caption{8D Volterra equation. From left to right: predicted solutions, exact solutions and the corresponding absolute errors.}
	\label{fig:volterra_8D_singular_kernel_result}
\end{figure}

\begin{figure}[!ht]
	\centering
	\includegraphics[scale=0.4]{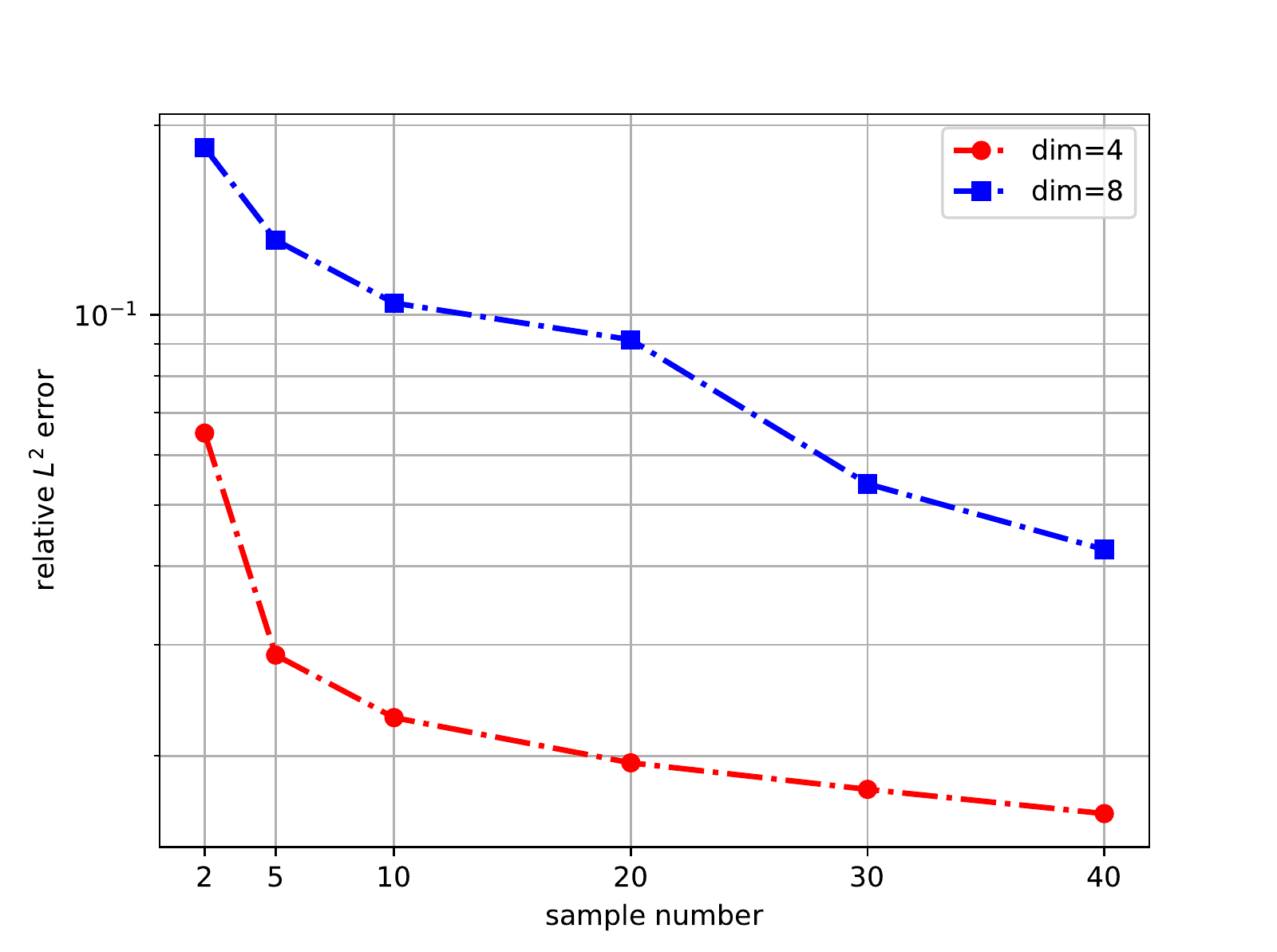}
	\caption{High dimensional volterra singular kernel problem. Relative $L^2$ error for different $N_s$.}
	\label{fig:volterra_high_dim_singualr_convergence}
\end{figure}

In this example, we take $d=3,7$ and 10 Monte Carlo points to approximate integral terms, i.e., $N_s=10$. The training set consists of two parts: 10000 collocation points randomly sampled in the equation domain and 1000 boundary points randomly sampled on each boundary. Then we train the MC-Nonlocal PINNs using LBFGS optimizer 40000 iterations. We select two different planes to elucidate the solution of our method. The final results are shown in Figure \ref{fig:volterra_4D_singular_kernel_result}, \ref{fig:volterra_8D_singular_kernel_result}. It is observed that the predictions achieve an excellent agreement with the corresponding ground truths for both $d=3$ and $d=7$. The corresponding relative $L^2$ errors for different sample numbers ($N_s$) are shown in Figure \ref{fig:volterra_high_dim_singualr_convergence}.

\subsection{Hypersingular integral equations}
\subsubsection{1D example}
\begin{figure}[H]
	\centering
	\includegraphics[scale=0.4]{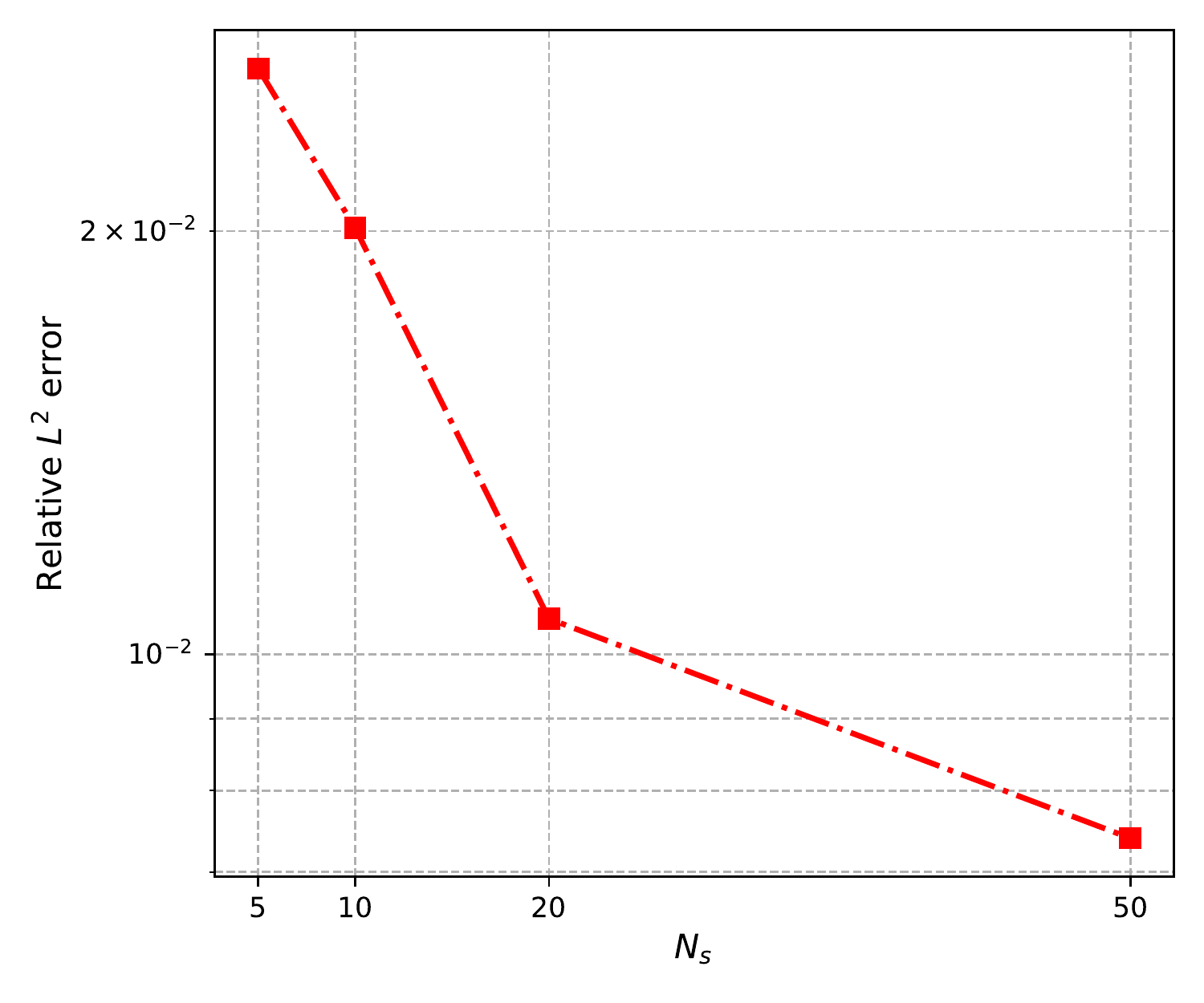}
	\caption{1D hypersingular integral equation. Relative $L^2$ errors for different $N_s$.}
	\label{fig:hypersingularconvergence}
\end{figure}

Consider the following hypersingular integral equation
\begin{equation}
	\int_0^1\hpsngAbs \frac{u(x)}{(x-s)^2}\mathrm{d}x= -\frac{1}{2} + 3s + (3s^2 - 2s)\ln \frac{1-s}{s},\quad \forall s\in(0,1),
\end{equation}
with boundary condition $u(0)=u(1)=0.$ The exact solution is $$u(x)=x^2(x-1).$$
We approximate the hypersingular integral via
\begin{equation}
	\frac{1}{N_s}\sum_{i=1}^{N_s}\left[(1-t_i)\frac{\partial ^2 u_{NN}}{\partial x^2}\big(t_ix_i + (1-t_i)s;\theta\big)\right] -u_{NN}(s;\theta)\bigg(\frac{1}{1-s} + \frac{1}{s}\bigg) + \frac{\partial u_{NN}}{\partial x}(s;\theta)\ln\frac{1-s}{s},
	\label{hyper_eqn}
\end{equation}
where $(t_i,x_i)\sim U[0,1]^2$ and $N_s$ is the number of discrete integration points. Note that Eq. (\ref{hyper_eqn}) holds for every $s\in(0,1)$, we use 100 uniformly distributed training points for $s$. We train the MC-Nonlocal-PINNs using the Adam optimizer for 2000 iterations, and relative $L^2$ errors between the exact and the predicted solutions for different $N_s$ are shown in Figure \ref{fig:hypersingularconvergence}.

\subsubsection{2D example}
\begin{figure}[H]
	\centering
	\includegraphics[scale=0.4]{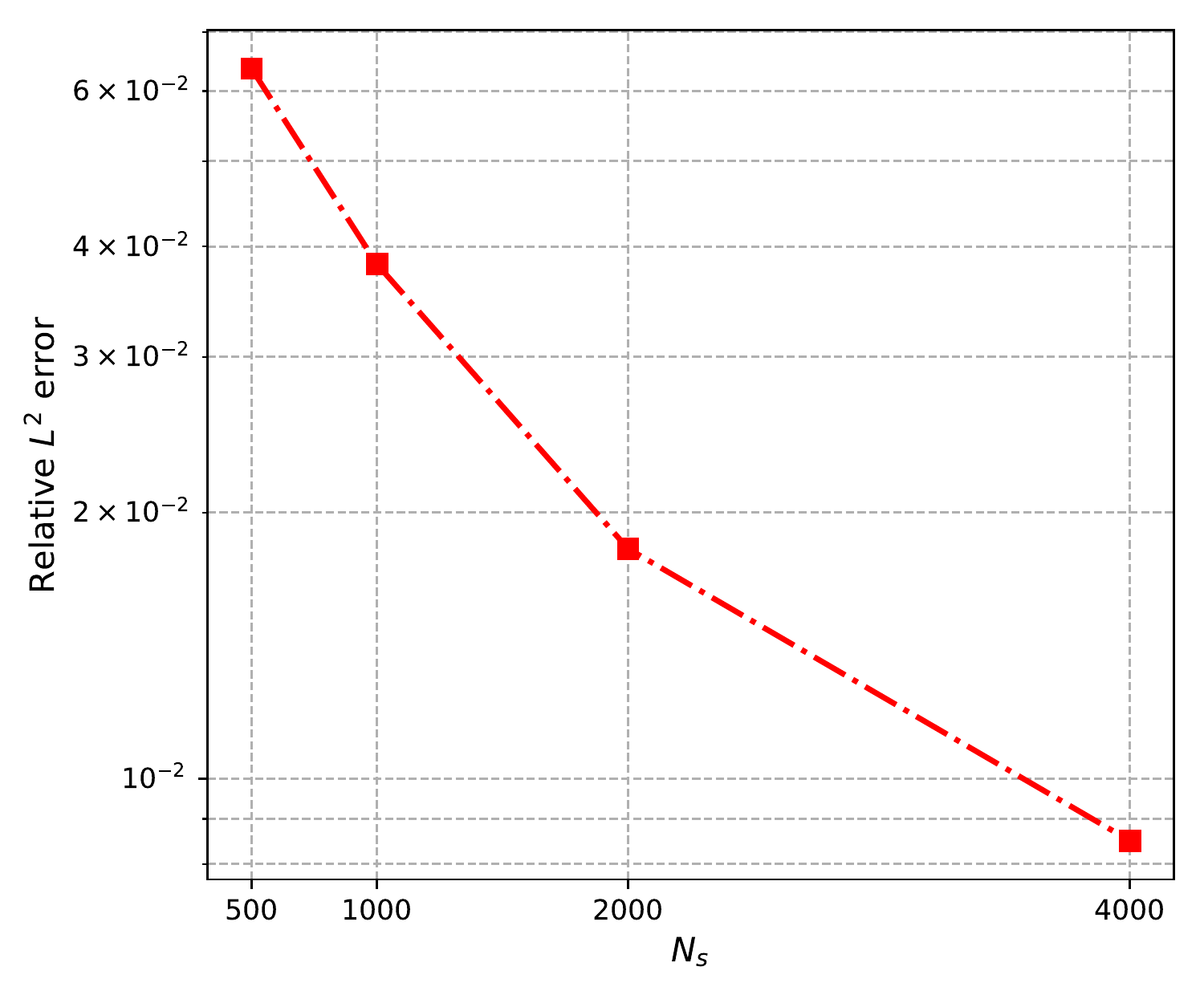}
	\caption{2D hypersingular integral equation. Relative $L^2$ errors for different $N_s$.}
	\label{fig:two_dim_hypersingularconvergence}
\end{figure}

Consider the following PDE involving hypersingular integral
\begin{equation}
	\begin{aligned}
		&-\Delta u + {\hiint}_\Omega \frac{u(x_1, x_2)}{r^3}dx_1dx_2 = f(x_1, x_2), & x\in \Omega,\\
		&u = g,& x\in \partial \Omega,
	\end{aligned}
\end{equation}
where $\Omega=\{(x_1,x_2)|x_1^2+x_2^2\leq 1\}$. The exact solution is given by
$$u(x_1,x_2) = \sin(\pi x_1)\sin(\pi x_2) + \exp(x_1+2x_2),$$
$f$ and $g$ can be calculated via classical numerical approach.
We approximate the above hypersingular integral via
\begin{equation*}
	{\hiint}_\Omega \frac{u(x_1,x_2)}{r^3}dx_1dx_2 \approx  \frac{2\pi}{N_s}\sum_{i=1}^{N_s}\left[-\tilde{u}_{NN}(0,\nu_i;\theta) + (1-t_i) \frac{\partial ^2 \tilde{u}_{NN}}{\partial r^2}(t_ir_i, \nu_i;\theta)
	\right],
\end{equation*}
where $\nu_i\sim \mc{U}[0,2\pi], t_i,r_i\sim \mc{U}[0,1]$ and $\tilde{u}_{NN}(r,\nu) = u_{NN}(r\cos\nu, r\sin\nu)$. And $\theta$ is learnable parameters.
We train the MC-Nonlocal-PINNs using the Adam optimizer for 40000 iterations. The final result is shown in Figure \ref{fig:two_dim_hypersingularconvergence}.

\subsection{Nonlocal PDEs}
\subsubsection{1D example}
Consider a one-dimensional nonlocal problem $-\mc{L}_{\delta}u=f_\delta$ on $(0,1)$ and the nonlocal operator is given by
$$\mathcal{L}_\delta u = 2\int _{-\delta}^{\delta}\gamma_\delta(s)(u(x+s)-u(x))\mathrm{d}s.$$
We mainly consider two cases: bound and singular kernel functions.

\begin{figure}[H]
	\centering
	\includegraphics[width=0.7\linewidth]{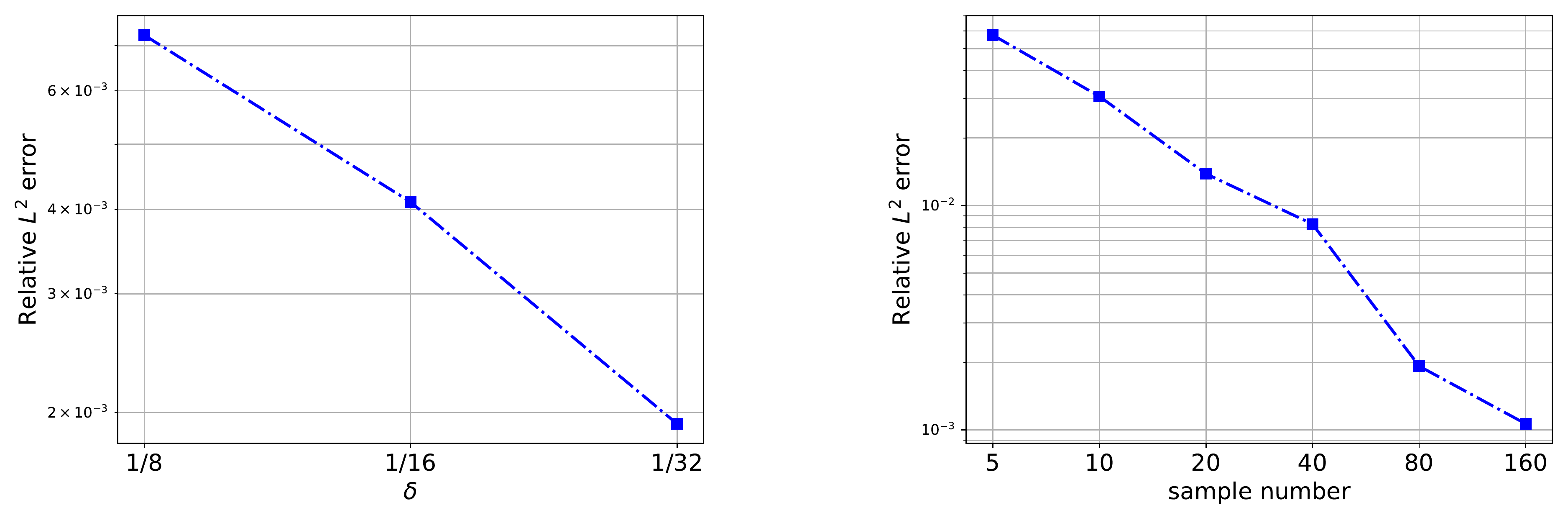}
	\caption{1D bounded kernel nonlocal problem. Relative $L^2$ error with different $\delta$ and $N_s$. Left : fixed sample number $N_s=80$. Right : fixed nonlocal radius $\delta=1/32$.}
	\label{fig:nonlocal1dbound}
\end{figure}

Case 1 (bounded kernel). A special kernel is chosen to be $\gamma_\delta(s)=\delta^{-2}|s|^{-1}$ in our numerical examples \cite{tian2013analysis}.
The exact solution is given by
$$u(x)=x^2(1-x^2),\quad f(x)= 12x^2-2+\delta^2.$$
And volume constraint is:
$$u(x)=x^2(1-x^2)\quad x\in (-\delta ,0)\cup (1,1+\delta).$$
For all $x\in(0,1)$, we can rewrite $\mc{L}_\delta$ as follows:
\begin{equation*}
	\begin{aligned}
		\mc{L}_\delta u(x) &= 2\int _0^{\delta}\delta^{-2} \frac{u(x+s)+u(x-s)-2u(x)}{s}\mathrm{d}s \approx \frac{1}{N_s}\sum_{i=1}^{N_s}\frac{u(x+\delta r_i) + u(x-\delta r_i) - 2u(x)}{r_i},
	\end{aligned}
\end{equation*}
where $r_i \sim U[0,1]$. The final result is shown in Figure \ref{fig:nonlocal1dbound}. We first fix sample number $N_s$ to be 80 and can observe that as the nonlocal radius $\delta$ decreases, the gap between our predicted solution and the reference becomes narrower. With fixed nonlocal radius $\delta=1/32$, relative $L^2$ error decreases as sample number ($N_s$) increases.

Case 2 (singular kernel). Another special kernel is chosen to be $\gamma(s)=\frac{1}{4}\delta^{-1/2}|s|^{-5/2}$ \cite{tian2015nonconforming}. Our benchmark problem is chosen to have $u(x)=-x^2(1-x)^2$ as the exact solution. The corresponding right-hand side $f(x)=12x^2-12x+2+\frac{2}{5}\delta^2$, and
$$u(x)=-x^2(1-x)^2\quad x\in (-\delta,0)\cup (1,1+\delta).$$

\begin{figure}[H]
	\centering
	\includegraphics[width=0.7\linewidth]{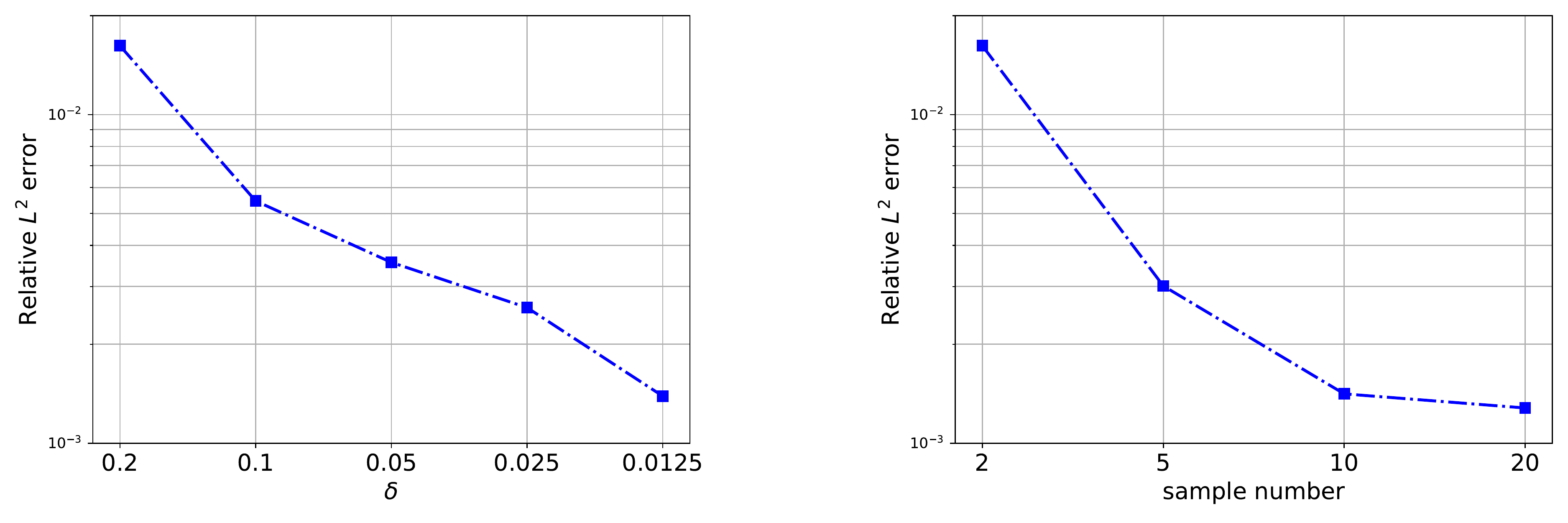}
	\caption{1D singular kernel nonlocal problem. Relative $L^2$ error with different $\delta$ and $N_s$. Left : fixed sample number $N_s=2$. Right : fixed nonlocal radius $\delta=0.2$.}
	\label{fig:nonlocal1dsingular}
\end{figure}

Similarly, we rewrite $\mathcal{L}_{\delta}$ as follows:
\begin{equation*}
	\begin{aligned}
		\mc{L}_{\delta}u(x)&=\frac{1}{2\sqrt{\delta}}\int_0^\delta
		s^{-1/2} \frac{u(x+s)-2u(x)+u(x-s)}{s^2}\mathrm{d}s\\
		& \approx \frac{1}{N_s}\sum_{i=1}^{N_s}\frac{u(x+\delta r_i) - 2 u(x) + u(x-\delta r_i)}{\delta^2 r_i^2},
	\end{aligned}
\end{equation*}	
where $r_i\sim \mbox{Beta}(0.5, 1)$. The final result is shown in Figure \ref{fig:nonlocal1dsingular}. We first fix sample number $N_s$ to 2 and can observe that as the nonlocal radius $\delta$ decreases, the gap between our predicted solution and the reference becomes narrower. With fixed nonlocal radius $\delta=0.2$, relative $L^2$ errors decrease as sample number ($N_s$) increases.

\subsubsection{High dimensional example}
\begin{figure}[H]
	\centering
	\includegraphics[scale=0.4]{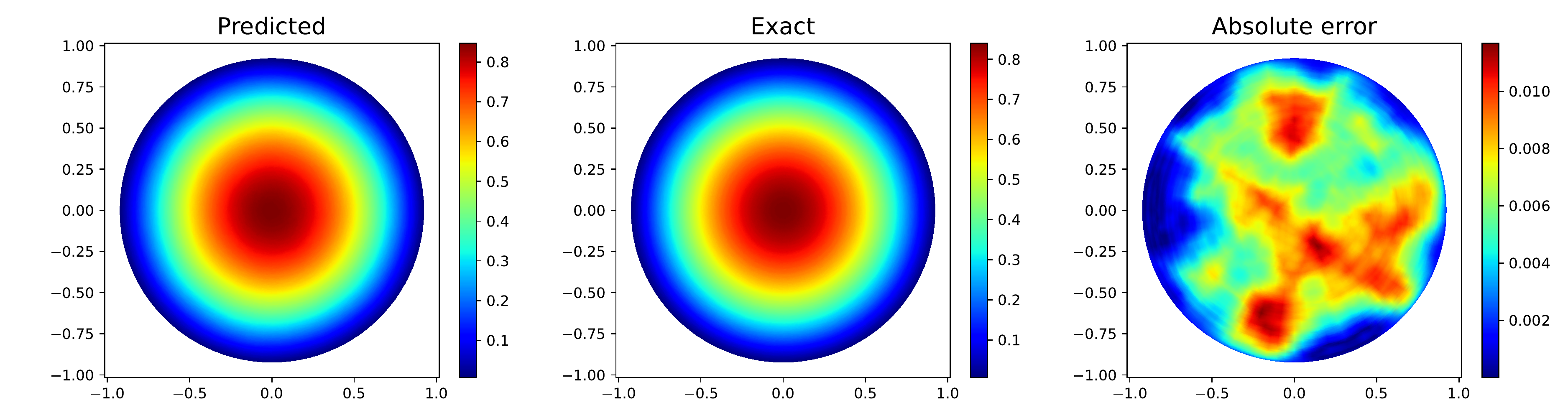}
	\caption{High dimensional nonlocal problem. From left to right: predicted solution, exact solution and absolute error.}
	\label{fig:highdimensionnonlocal}
\end{figure}

In this example, we consider a four dimensional nonlocal problem with Dirichlet boundary condition:

\begin{equation}
	\left\{	\begin{aligned}
		&	\int_{B(x,\delta)}\frac{u(x)-u(y)}{\Vert x-y\Vert_2^{d+\alpha}}dy = f(x),& \mbox{in }\Omega,\\
		&u(x)=g(x),& \mbox{in } \Omega_\delta.
	\end{aligned}
	\right.
\end{equation}
The fabricated solution is
\begin{equation}
	u(x) = (1-\Vert x\Vert_2^2)^{\alpha/2},\quad x\in \Omega = \mathbb{B}_1^{4}=\left\{x\mid \Vert x\Vert_2\leq 1, x\in \mathbb{R}^4\right\},
\end{equation}
the corresponding force term $f(x)$ can be calculated via classical numerical approach. Here we take homogeneous nonlocal boundary condition, that is, $g(x)=0.$ We set $\delta=0.2, \alpha=0.5$. We use ReLU as activation function and approximate $u(x)$ with $u_{NN}(x)=\mbox{ReLU}(1-\Vert x\Vert ^2)\tilde{u}_{NN}(x)$ to exactly satisfy non-zero volume boundary condition, and train the MC-Nonlocal-PINNs using the Adam optimizer for 40000 iterations with batch size 128. To illustrate final result, we select one plane $[x_1, x_2, 0.2, 0.2]$. The Figure \ref{fig:highdimensionnonlocal} shows the comparison between predicted and exact solutions. The relative $L^2$ error is 6.85e-3, which is sufficiently low for the high dimensional problem.

\section{Summary}

We have proposed MC-Nonlocal-PINNs for solving general nonlocal models such as volterra-type, hypersingular integral equations and nonlocal PDEs. Our MC-Nonlocal-PINNs handle the nonlocal operators in a Monte Carlo way, resulting in a very stable approach for high dimensional problems. Applications to hypersingular integral equations, high dimensional Volterra type integral equations and nonlocal PDEs demonstrate the effectiveness of our approach.

Despite the encouraging results presented here, some integral equations  still require further investigation such as highly oscillatory kernels arsing in electromagnetics, integration on manifolds and complex domains. We believe that addressing these problems will provide a better understanding of the MC-Nonlocal-PINNs approach.

%%%% Acknowledgments %%%%%%%%
\section*{Acknowledgments}
This study was sponsored by the National Natural Science Foundation of China (NSFC: 11971259).

%%%% Bibliography  %%%%%%%%%%
%\nobibliography{ref.bib}

\end{document}